\documentclass[reqno,12pt]{amsart}
\usepackage{amsmath,amsthm,amssymb,amsfonts,amscd,graphicx}
\input xy
\xyoption{all}
\theoremstyle{plain} 
	\newtheorem{thm}{Theorem}[section]
	\newtheorem*{thm*}{Theorem}
	\newtheorem{cor}[thm]{Corollary}
	\newtheorem{lem}[thm]{Lemma}
	\newtheorem{sublem}[thm]{Sub-Lemma}
	\newtheorem{prop}[thm]{Proposition}
	\newtheorem*{prop*}{Proposition}
	\newtheorem{conj}[thm]{Conjecture}
	\newtheorem*{conj*}{Conjecture}
	
\theoremstyle{definition}
	\newtheorem{defn}[thm]{Definition}

\theoremstyle{remark}
	\newtheorem{rem}[thm]{Remark}
	
	\newtheorem*{pf}{Proof}
\numberwithin{equation}{section}

\def\CC{{\mathbb C}}

\def\HH{{\mathbb H}}
\def\LL{{\mathbb L}}
\def\NN{{\mathbb N}}
\def\PP{{\mathbb P}}
\def\QQ{{\mathbb Q}}
\def\RR{{\mathbb R}}
\def\TT{{\mathbb T}}
\def\ZZ{{\mathbb Z}}
\def\A{{\mathcal A}}

\def\C{{\mathcal C}}
\def\D{{\mathcal D}}
\def\E{{\mathcal E}}

\def\I{{\mathcal I}}

\def\L{{\mathcal L}}

\def\R{{\mathcal R}}

\def\S{{\mathcal S}}

\def\g{{\mathfrak g}}

\begin{document}
\title{On Weyl Groups and Artin Groups Associated to Orbifold Projective Lines}
\date{\today}
\author{Yuuki Shiraishi}
\address{Department of Mathematics, Graduate School of Science, Osaka University, 
Toyonaka Osaka, 560-0043, Japan}
\email{sm5021sy@ecs.cmc.osaka-u.ac.jp}
\author{Atsushi Takahashi}
\address{Department of Mathematics, Graduate School of Science, Osaka University, 
Toyonaka Osaka, 560-0043, Japan}
\email{takahashi@math.sci.osaka-u.ac.jp}
\author{Kentaro Wada}
\address{Department of Mathematics, Graduate School of Science, Osaka University, 
Toyonaka Osaka, 560-0043, Japan}
\email{k-wada@cr.math.sci.osaka-u.ac.jp}
\begin{abstract}
We associate a generalized root system in the sense of Kyoji Saito to an orbifold projective line via the derived category of finite dimensional representations of a certain bound quiver algebra.
We generalize results by Saito--Takebayshi and Yamada for elliptic Weyl groups and elliptic Artin groups 
to the Weyl groups and the fundamental groups of the regular orbit spaces associated to the generalized root systems.
Moreover we study the relation between this fundamental group and a certain subgroup of the autoequivalence group of a triangulated
subcategory of the derived category of $2$-Calabi--Yau completion of the bound quiver algebra.   
\end{abstract}
\maketitle
\section{Introduction}
The purpose of the present paper is to associate a generalized root system in the sense of Kyoji Saito to an orbifold projective line and to understand some properties of the Weyl group of the generalized root system.
In particular, we want to know the description of this Weyl group as the group with generators and relations
and the connection between the Weyl group and the fundamental group of the regular orbit space associated to the generalized root system.
Saito--Takebayashi \cite{st:1} introduces the generalized Coxeter relations in order to describe by generators and relations 
an elliptic Weyl group, 
the Weyl group of an elliptic root system which is an invariant of the mirror dual object to an elliptic orbifold projective line.
It is known that the base space of a universal unfolding of a simple elliptic singularity 
is described as the orbit space associated to the corresponding elliptic root system and, in particular, the complement of the discriminant of the base space is homeomorphic to the regular orbit space (cf. Chapter 3 in \cite{looijenga}). 
Yamada \cite{yamada} gives the generalized Coxeter relations to the fundamental group of this regular orbit space
and shows that there is a natural surjective map from the fundamental group of it to the Weyl group.
On the other hand, one of the most important theme in mirror symmetry is to understand the relations 
among complex geometry, symplectic geometry and representation theory via the isomorphisms of 
Frobenius structures constructed from these mathematical areas.
The Frobenius structure is defined by Dubrovin \cite{d:1} in order to obtain the geometric description of the WDVV equations 
on the topological field theory. On the other hand, the Frobenius structure was first discovered by K. Saito as the flat structure on the base space of a universal unfolding of an isolated hypersurface singularity in his study of a primitive form \cite{saito} 
and on a complexified quotient space of a real Euclid space by the finite reflection group on it \cite{saito:1,ssy}.
The notion of the generalized root system is introduced in \cite{saito}
in order to understand algebraically vanishing cycles in the Milnor fiber and 
the period mapping of a primitive form.
From the view point of mirror symmetry, the generalized root system 
should exist on other mathematical areas.
Moreover the generalized root system introduced in other areas
should play an important role to understand them
as vanishing cycles do in the singularity theory.
Let us give the plan of this paper and then explain more precisely our motivation. 
Let $r\ge 3$ be a positive integer, $A = (a_{1},\ldots,a_{r})$ an  $r$-tuple of positive integers greater than one and 
$\Lambda=(\lambda_1,\ldots,\lambda_r)$ an $r$-tuple of pairwise distinct elements of $\PP^1(k)$ normalized such that 
$\lambda_1=\infty$, $\lambda_2=0$ and $\lambda_3=1$. 
We consider an orbifold projective line $\PP^{1}_{A,\Lambda}$ of type $(A,\Lambda)$ (see Definition~\ref{def:2.22}), 
which is the smooth proper Deligne--Mumford stack 
whose coarse moduli space is $\PP^{1}(k)$ with the orbifold points $\lambda_{i}$ of orders $a_{i}$ for $i=1,\dots, r$. 
Starting from an orbifold projective line $\PP^{1}_{A,\Lambda}$ and passing through the bounded derived category $\D^b{\rm coh}(\PP^1_{A,\Lambda})$ of coherent sheaves on $\PP^1_{A,\Lambda}$, 
we shall associate to it a generalized root system $\widetilde{R}_A$, which is independent from $\Lambda$, based on Proposition~\ref{prop:1} which says that a generalized root system is an invariant for a triangulated category satisfying some conditions 
which is expected to hold for triangulated categories of a good algebro-geometric origin \cite{BP, CB1, kuleshov-orlov, ringel, rudakov}.
One of key facts is the triangulated equivalence between $\D^b{\rm coh}(\PP^1_{A,\Lambda})$ and 
the bounded derived category $\D^{b}(k\widetilde{\TT}_{A,\Lambda})$ of finite dimensional representations of the bound quiver algebra $k\widetilde{\TT}_{A,\Lambda}$ (see Definition \ref{defn of widetilde{TT}_A}), which is a direct consequence of the result by Geigle--Lenzing \cite{gl:1}.
The reason why we consider the bound quiver $\widetilde{\TT}_{A,\Lambda}$ is that
we can obtain the generalized Coxeter relations introduced by Saito--Takebayashi \cite{st:1} and Yamada \cite{yamada}.
At this moment, we have not succeeded yet to obtain a generalization of Coxeter relations for Ringel's canonical algebra \cite{ringel:1} which is also derived equivalent to $\D^b{\rm coh}(\PP^1_{A,\Lambda})$. 
The results for elliptic Weyl groups in \cite{st:1} are naturally extended.
We describe the {\it cuspidal Weyl group}, the Weyl group $W(\widetilde{R}_A)$ of $\widetilde{R}_A$, as a semi-direct product of the Weyl group $W(R_A)$ of $R_A$ and the Grothendieck group $K_{0}(R_A)$ 
(see Theorem~\ref{relation W(R_A) and W(widetilde{R}_A)} and Definition~\ref{defn:cuspidal Weyl}), 
where $R_A$ is a generalized root system associated to the bounded derived category $\D^{b}(k\TT_A)$ of finite dimensional representations of the path algebra $k\TT_A$ for the subquiver $\TT_A$ of $\widetilde{\TT}_A$ (see Definition \ref{defn:2.19}). 
Recall that, for elliptic cases, the hyperbolic extension of $W(\widetilde{R}_A)$ is isomorphic to $W(R_A)\ltimes K_{0}(R_A)$.
We check in Theorem \ref{main} the generalized Coxeter relations in \cite{st:1}, a generalization of Coxeter relations for the bound components in Coxeter--Dynkin diagram,
is also valid for the Weyl group $W(\widetilde{R}_A)$.
The results for elliptic Artin groups in \cite{yamada} are also naturally extended.
We define a group $G(\widetilde{T}_A)$ (see Definition~\ref{def G'}) for 
the Coxeter--Dynkin diagram $\widetilde{T}_A$, which is called the {\it cuspidal Artin group} of type $A$.
Since the group $W(R_A)\ltimes K_{0}(R_A)$ naturally acts on the complexified Tits cone $\E(R_A)$ associated to ${R}_A$ in a properly discontinuous way and its action is free on the complement set $\E(R)^{reg}$ of the reflection hyperplanes, we can consider the regular orbit space $\E(R)^{reg}/(W(R_A)\ltimes K_{0}(R_A))$.
We show in Theorem \ref{G isom G'} that $G(\widetilde{T}_A)$ is 
isomorphic to the fundamental group $G(\widetilde{R}_A)$ of 
the regular orbit space $\E(R)^{reg}/(W(R_A)\ltimes K_{0}(R_A))$. 
From now on we shall return to our motivation in mirror symmetry. 
In particular, the generalized root system $\widetilde{R}_A$ will be related to the space of stability conditions for the derived category. 
The second named author expects the following conjecture.
For simplicity, we restrict ourselves to the non-elliptic case:
\begin{conj}\label{conj:1.1}
Set $\chi_A:=2+\displaystyle\sum_{i=1}^r (1/a_i-1)$ and assume that $\chi_A\ne 0$.
\begin{enumerate}
\item
The complex manifold 
\begin{equation}
M_{A}:= 
\begin{cases}
(\E(R_A)\times \CC)/W(R_A)\ltimes K_{0}(R_A) \quad \text{if} \quad \chi_A>0,\\
(\E(R_A)\times \HH)/W(R_A)\ltimes K_{0}(R_A)\quad \text{if} \quad \chi_A<0
\end{cases}
\end{equation}
should be isomorphic to the space of stability conditions ${\rm Stab}(\D^{b}{\rm coh}(\PP^{1}_{A,\Lambda}))$ for $\D^{b}{\rm coh}(\PP^{1}_{A,\Lambda})$.
Here the group $W(R_A)\ltimes K_{0}(R_A)$ acts on $\E(R_A)$ as reflections and translations, and on $\CC$ or $\HH$ trivially.
\item 
There should exist a Frobenius structure on $M_A$ which 
is locally isomorphic to the one constructed from the Gromov--Witten theory of $\PP^{1}_{A,\Lambda}$. 
\item Set the complex manifold
\begin{equation}
M^{reg}_{A}:=
\begin{cases}
(\E(R_A)^{reg} \times \CC)/W(R_A)\ltimes K_{0}(R_A) \quad \text{if} \quad \chi_A>0,\\
(\E(R_A)^{reg} \times \HH)/W(R_A)\ltimes K_{0}(R_A) \quad \text{if} \quad \chi_A<0.
\end{cases}
\end{equation}
Then the universal covering $\widetilde{M}^{reg}_A$ of $M^{reg}_A$ should be the space of stability conditions 
${\rm Stab}(\check{\D}_{A,\Lambda})$ for $\check{\D}_{A,\Lambda}$ where $\check{\D}_{A,\Lambda}$ is 
the smallest full triangulated subcategory of the derived category of the $2$-Calabi--Yau completion $\Pi_{2}(k\widetilde{\TT}_{A,\Lambda})$ 
of $k\widetilde{\TT}_{A,\Lambda}$ containing $k\widetilde{\TT}_{A,\Lambda}$, closed under isomorphisms and taking direct summand.
\item 
Assume that $r=3$. By \cite{ist:2, ist:3} and {\rm (ii)},
the Frobenius structure on $M_A$ in {\rm (ii)} can be identified with the one on 
the $\ZZ$-covering of the base space of a universal unfolding of a cusp polynomial $f_A:=x_{1}^{a_1}+x_{2}^{a_2}+x_{3}^{a_3}-s^{-1}_{\mu_A}x_1x_2x_3$ with 
the primitive form $\zeta^{(0)}:=\zeta_A$ obtained in \cite{ist:2} for $\chi_A>0$ and in \cite{ist:3} for $\chi_A<0$.
Under this identification,
the stability function 
\begin{equation}
Z: {\rm Stab}(\check{\D}_{A,\Lambda})\longrightarrow (K_{0}(\check{\D}_{A,\Lambda})\otimes_\ZZ \CC)^{*}
\end{equation}
should be given by the period mapping$:$ 
\begin{equation}
\int \check{\zeta}^{(-1)}: \widetilde{M}^{reg}_A
\longrightarrow (K_{0}(\check{\D}_{A,\Lambda})\otimes_\ZZ \CC)^{*},
\end{equation}
where the element $\check{\zeta}^{(-1)}$ is the formal Fourier--Laplace transform $($the Gelfand--Leray form$)$ of the primitive form $\zeta^{(-1)}$ twisted by $-1$ $($cf. Section 5 of \cite{S1202-Saito}$)$.
\end{enumerate}
\end{conj}
We shall explain some known results concerning this conjecture.
Looijenga shows in Chapter 3 of \cite{looijenga} that the base space of a universal unfolding of $f_A$ for $r=3$ and $\chi_A\le 0$ can be described as the orbit space in Conjecture \ref{conj:1.1} {\rm (i)} as complex manifolds.  
In particular, the complement of the discriminant of this universal unfolding of $f_A$ is homeomorphic to the regular orbit space in Conjecture \ref{conj:1.1} {\rm (iii)} 
as complex manifolds. 
If $\chi_A>0$, Dubrovin--Zhang \cite{dz:1} constructed a Frobenius structure on $M_A$. 
It is shown in Ishibashi--Shiraishi--Takahashi \cite{ist:1, ist:2} and Shiraishi--Takahashi \cite{ist:3} that
this Frobenius structure is isomorphic to the one constructed from the Gromov--Witten theory of $\PP^{1}_{A,\Lambda}$
and the one constructed from the universal unfolding of the cusp polynomial $f_A$ with the primitive form $\zeta_A$. 
The group $G(\widetilde{R}_A)$ is related to a subgroup of the autoequivalence group ${\rm Auteq}(\check{\D}_{A,\Lambda})$.
Let ${\rm Br}(\check{\D}_{A,\Lambda})$ be the subgroup of the autoequivalence group
${\rm Auteq}(\check{\D}_{A,\Lambda})$ generated by spherical twist functors, defined in Seidel--Thomas \cite{Seidel-Thomas}, 
of simple $k\widetilde{\TT}_{A,\Lambda}$-modules regarded as $\Pi_{2}(k\widetilde{\TT}_{A,\Lambda})$-modules. 
Then we have the natural surjective group homomorphism from $G(\widetilde{R}_A)$ to ${\rm Br}(\check{\D}_{A,\Lambda})$ in Theorem \ref{isom G and Br}. 
Under Conjecture \ref{conj:1.1} {\rm (iii)} and {\rm (iv)}, this surjectivity can be viewed as the property of the covering map and the fundamental group
since ${\rm Br}(\check{\D}_{A, \Lambda})$ is regarded as the group of covering transformations for ${\rm Stab}(\check{\D}_{A,\Lambda})$.
Similar to the results by Bridgeland for K3 surfaces in \cite{Bri:1} and Kleinian singularities in \cite{Bri:2}, we expect the following conjecture:
\begin{conj}\label{conj:1.2}
The group homomorphism $G(\widetilde{T}_A)\twoheadrightarrow {\rm Br}(\check{\D}_{A,\Lambda})$ in Theorem \ref{isom G and Br} 
should also be injective, and hence isomorphism.
In other words, the space of stability condition ${\rm Stab}(\check{\D}_{A,\Lambda})$ should be simply connected. 
\end{conj}
Similar known results for
the injectivity of the group homomorphism in Conjecture \ref{conj:1.2} are obtained 
by Brav--Thomas \cite{bt:1}, Ishii--Ueda--Uehara \cite{iuu:1} and Seidel--Thomas \cite{Seidel-Thomas}.
The above conjecture is the further theme to be worked on with Conjecture \ref{conj:1.1}.
\bigskip
\noindent
{\it Acknowledgement}\\
\indent
The second named author is supported by JSPS
KAKENHI Grant Number 24684005.
\section{Notations and terminologies}
Throughout this paper, $k$ denotes an algebraically closed field of characteristic zero.
\subsection{Generalized root systems}
In this subsection, we recall the definition of the simply laced generalized root system
introduced by K. Saito \cite{saito}.
\begin{defn}
A {\it simply-laced generalized root system} $R$ consists of  
\begin{itemize}
\item
a free $\ZZ$-module $K_0(R)$ of finite rank ($=:\mu$) called the {\it root lattice},
\item
a symmetric bi-linear form $I_R:K_0(R)\times K_0(R)\longrightarrow\ZZ$,
\item
a subset $\Delta_{re}(R)$ of $K_0(R)$ called the {\it set of real roots} such that
\begin{enumerate}
\item
$K_0(R)=\ZZ\Delta_{re}(R)$,
\item
for all $\alpha\in\Delta_{re}(R)$, $I(\alpha,\alpha)=2$,
\item 
for all $\alpha\in\Delta_{re}(R)$, the element $r_\alpha$ of ${\rm Aut}(K_0(R),I_R)$, the group of automorphisms of $K_0(R)$ respecting $I_R$, defined by
\begin{equation}
r_\alpha(\lambda):=\lambda-I_R(\lambda,\alpha)\alpha,\quad \lambda\in K_0(R),
\end{equation}
makes $\Delta_{re}(R)$ invariant, namely, $r_\alpha(\Delta_{re}(R))=\Delta_{re}(R)$,
\item 
there exists a subset $B=\{\alpha_1,\dots, \alpha_\mu\}$ of $\Delta_{re}(R)$ called a {\it root basis} of $R$ which satisfies 
$K_0(R)=\displaystyle\bigoplus_{i=1}^\mu \ZZ\alpha_i$, $W(R)=\langle r_{\alpha_1},\dots,r_{\alpha_\mu}\rangle$ and $\Delta_{re}(R)=W(R)B$ where $W(R)$ is the {\it Weyl group} of $R$ defined by
\begin{equation}
W(R):=\langle r_{\alpha}~\vert~\alpha\in\Delta_{re}(R)\rangle\subset {\rm Aut}(K_0(R),I_R), 
\end{equation}
\end{enumerate}
\item
an element $c_R$ of $W(R)$ called the {\it Coxeter transformation} which has the presentation $c_R=r_{\alpha_1}\cdots r_{\alpha_\mu}$ with respect to a root basis $B$. 
\end{itemize}
An element of $\Delta_{re}(R)$ is called a {\it real root} and an element of $B$ is called a {\it real simple root}.  
For a real simple root $\alpha\in B$, the reflection $r_{\alpha}$ is called a {\it simple reflection}.
\end{defn}
One can define a notion of isomorphism of simply-laced generalized root systems in the obvious way. 
\begin{defn}
Let $R=(K_0(R),I_R,\Delta_{re}(R),c_R)$ be a simply-laced generalized root system with a root basis $B=\{\alpha_1,\dots, \alpha_\mu\}$ of $R$. 
The {\it Coxeter--Dynkin diagram} $\Gamma_B$ is a finite graph defined as follows$:$ 
\begin{itemize}
\item
the set of vertices is $B=\{\alpha_1,\dots, \alpha_\mu\}$,
\item
the edge between vertices $\alpha_i$ and $\alpha_j$ of $\Gamma_B$ is given by the following rule$:$
\begin{subequations}
\begin{eqnarray}
\xymatrix{ \circ_{\alpha_i}  & \circ_{\alpha_j}} & \quad\text{if}\quad I_R(\alpha_i,\alpha_j)=0,\\
\xymatrix{ \circ_{\alpha_i}\ar@{-}[r]  & \circ_{\alpha_j}} & \quad\text{if}\quad I_R(\alpha_i,\alpha_j)=-1,\\
\xymatrix{ \circ_{\alpha_i}\ar@{-}[r]_{t}  & \circ_{\alpha_j}} & \quad\text{if}\quad I_R(\alpha_i,\alpha_j)=-t,\ (t\ge 2),\\
\xymatrix{ \circ_{\alpha_i}\ar@{.}[r]  & \circ_{\alpha_j}} & \quad\text{if}\quad I_R(\alpha_i,\alpha_j)=+1,\\
\xymatrix{ \circ_{\alpha_i}\ar@2{.}[r]  & \circ_{\alpha_j}} & \quad\text{if}\quad I_R(\alpha_i,\alpha_j)=+2,\\
\xymatrix{ \circ_{\alpha_i}\ar@{.}[r]_{t}  & \circ_{\alpha_j}} & \quad\text{if}\quad I_R(\alpha_i,\alpha_j)=+t,\ (t\ge 3).
\end{eqnarray}
\end{subequations}
\end{itemize}
\end{defn}
\subsection{Generalized root systems from triangulated categories}

In this subsection, we deduce a simply laced generalized root system from a certain algebraic triangulated category
which satisfies plausible conditions.
\begin{defn}
Let $\D$ be a $k$-linear triangulated category with the translation functor $[1]$.
Consider a free abelian group $F$ with generators $\{[X]~|~X\in\D\}$ and a subgroup $F_0$ of $F$ generated by $[X]-[Y]+[Z]$ for all exact triangles $X\longrightarrow Y\longrightarrow Z\longrightarrow X[1]$ in $\D$.
The {\it Grothendieck group} $K_0(\D)$ of $\D$ is a quotient group $F/F_0$.
\end{defn}
Any triangulated category of our interest in this paper is equipped with an enhancement. 
We briefly recall some terminologies.
\begin{defn}[Keller~\cite{kellerdg}]
Let $\D$ be a $k$-linear triangulated category.
We say that $\D$ is {\it algebraic} if it is equivalent as a triangulated category to the stable category of some $k$-linear 
Frobenius category. 
\end{defn}
It is important to note that for an algebraic $k$-linear triangulated category $\D$, 
we have functorial cones and $\RR{\rm Hom}$-complexes once 
we fix an enhancement, a differential graded category which yields $\D$ 
(see Theorem~3.8 in \cite{kellerICM} for precise statements).
\begin{defn}
Let $\D$ be an algebraic $k$-linear triangulated category with the translation functor $[1]$ with a fixed enhancement.
\begin{enumerate}
\item
For $X,Y\in\D$, denote by $\RR{\rm Hom}_\D^{\bullet}(X,Y)\in \D(k)$ the $\RR{\rm Hom}$-complex 
such that ${\rm Hom}_\D(X,Y[p])=H^p(\RR{\rm Hom}_\D^{\bullet}(X,Y))$ for all $p\in\ZZ$, where 
$\D(k)$ is the derived category of complexes of $k$-modules.
\item
A $k$-linear triangulated category $\D$ is said to be {\it of finite type} if the total dimension of the graded $k$-module
${\rm Hom}_\D^{\bullet}(X,Y):=\displaystyle\bigoplus_{p\in\ZZ}{\rm Hom}_\D(X,Y[p])[-p]$ is finite for all $X,Y\in\D$.
\end{enumerate}
\end{defn}
\begin{defn}
Let $\D$ be an algebraic $k$-linear triangulated category of finite type with a fixed enhancement.
\begin{enumerate}
\item
An object $E$ in $\D$ is called an {\it exceptional object} (or is called {\it exceptional}) if 
$\RR{\rm Hom}_\D^{\bullet}(E,E)\cong k\cdot {\rm id}_E$ in $\D(k)$.
\item 
An {\it exceptional collection} $\E=(E_1,\dots, E_n)$ in $\D$ is a finite ordered set of exceptional objects 
satisfying the condition that $\RR{\rm Hom}_\D^{\bullet}(E_i,E_j)\cong 0$ in $\D(k)$ for all $i>j$.
An exceptional collection consisting of two objects is an {\it exceptional pair}.
\item 
An exceptional collection $\E=(E_1,\dots, E_n)$ in $\D$ is said to be {\it isomorphic} to another 
exceptional collection $\E'=(E'_1,\dots, E'_n)$ in $\D$ if $E_i\cong E'_i$ in $\D$ for all $i=1,\dots, \mu$.
\item
An exceptional collection $\E=(E_1,\dots, E_n)$ in $\D$ is called 
a {\it strongly exceptional collection} if, for all $i,j=1,\dots, n$, the complex $\RR{\rm Hom}_\D^{\bullet}(E_i,E_j)$ is isomorphic in $\D(k)$ to a complex concentrated in degree zero, equivalently, we have ${\rm Hom}_{\D}(E_i,E_j[p])=0$ for $p\ne 0$.
\item 
An exceptional collection $\E$ in $\D$ is called {\it full}
if the smallest full triangulated subcategory of $\D$ containing all elements in $\E$
is equivalent to $\D$.
\item
For an exceptional pair $(X,Y)$, one has new exceptional pairs $(\L_{X}Y, Y)$ called the {\it left mutation} of $(X,Y)$ 
and $(Y, \R_{Y}X)$ called the {\it right mutation} of $(X,Y)$.
Here the object $\L_XY[1]$ is defined as the cone of the evaluation morphism $ev$
\begin{subequations}
\begin{equation}
\RR{\rm Hom}_\D^{\bullet}(X,Y)\otimes^\LL X \stackrel{ev}{\longrightarrow} Y,
\end{equation}
where $(-)\otimes^{\LL} X$ is the left adjoint of the functor $\RR{\rm Hom}_\D(X,-):\D\longrightarrow \D(k)$.
Similarly, the object $\R_Y X$ is defined as the cone of the evaluation morphism $ev^*$
\begin{equation}
X \stackrel{ev^*}{\longrightarrow} \RR{\rm Hom}_\D^{\bullet}(X,Y)^{*}\otimes^\LL Y.
\end{equation}
\end{subequations}
where $(-)^*$ denotes the duality ${\rm Hom}_k(-,k)$.
\end{enumerate}
\end{defn}
Here we recall the braid group action on the set of  isomorphism classes of full exceptional collections.
\begin{defn}
The Artin's {\it braid group} $B_{\mu}$ on $\mu$-stands is a group presented by the following generators and relations: 
\begin{description}
\item[{\bf Generators}] $\{b_i~|~i=1,\dots, \mu-1\}$
\item[{\bf Relations}] 
\begin{subequations}
\begin{equation}
b_{i}b_{j}=b_{j}b_{i}\quad \text{for}\quad|i-j|\ge 2,
\end{equation}
\begin{equation}
b_{i}b_{i+1}b_{i}=b_{i+1}b_{i}b_{i+1}\quad \text{for}\quad i=1,\dots, \mu-2.
\end{equation}
\end{subequations}
\end{description}
\end{defn}
Consider the group $G_\mu:=B_{\mu}\ltimes \ZZ^{\mu}$, the semi-direct product of the braid group $B_\mu$ and the 
free abelian group of rank $\mu$, defined by the group
homomorphism $B_{\mu}\rightarrow {\mathfrak S}_{\mu} \rightarrow {\rm Aut}_{\ZZ} \ZZ^{\mu}$,
where the first homomorphism is $b_{i}\mapsto (i, i+1)$ and the second one is
induced by the natural actions of the symmetric group ${\mathfrak S}_{\mu}$ on $\ZZ^{\mu}$.  
\begin{prop}[cf. Proposition~2.1 in \cite{BP}]
Let $\D$ be an algebraic $k$-linear triangulated category of finite type with a fixed enhancement..
The group $G_\mu$ acts on the set of isomorphism classes of full exceptional collections in $\D$ by mutations and translations$:$
\begin{subequations}
\begin{equation}
b_{i}(E_{1},\dots, E_{\mu}):=(E_{1}, \dots, E_{i-1}, E_{i+1}, \R_{E_{i+1}}E_i,E_{i+2},\dots, E_{\mu}),
\end{equation}
\begin{equation}
b^{-1}_{i}(E_{1},\dots, E_{\mu}):=(E_{1}, \dots, E_{i-1}, \L_{E_{i}}E_{i+1}, E_{i}, E_{i+2},\dots, E_{\mu}),
\end{equation}
\begin{equation}
e_{i}(E_{1},\dots, E_{\mu}):=(E_{1}, \dots, E_{i-1}, E_{i}[1], E_{i+1},\dots, E_{\mu}),
\end{equation}
where we denote by $e_{i}$ the $i$-th standard basis of $\ZZ^{\mu}$.
\end{subequations}
\qed
\end{prop}
\begin{prop}\label{prop:1}
Let $\D$ be an algebraic $k$-linear triangulated category of finite type with the translation functor $[1]$ and a fixed enhancement.
Assume that $\D$ satisfies the following conditions$:$
\begin{enumerate}
\item
There exists a full strongly exceptional collection $\E=(E_1,\dots,E_\mu)$ in $\D$.
\item
The action of the group $G_\mu$ on the set of isomorphism classes of full exceptional collections in $\D$ is transitive.
\item
For any exceptional object $E'\in \D$, there exists a full exceptional collection $\E'$ in $\D$ such that $E'\in\E'$.
\end{enumerate}
Then the following quadruple
\begin{itemize}
\item
the Grothendieck group $K_0(\D)$ of $\D$,
\item
the Cartan form $I_\D:K_0(\D)\times K_0(\D)\longrightarrow \ZZ;$
\begin{equation}
I_\D([X],[Y]):=\chi_\D([X],[Y])+\chi_\D([Y],[X]),\quad X,Y\in\D,
\end{equation}
where $\chi_\D:K_0(\D)\times K_0(\D)\longrightarrow \ZZ$ is the Euler form defined by  
\begin{equation}
\chi_\D([X],[Y]):=\sum_{p\in \ZZ}(-1)^p\dim_k {\rm Hom}_{\D}(X,Y[p]),
\end{equation}
\item
the subset $\Delta_{re}(\D)$ of $K_0(\D)$ defines by  
\begin{equation}
\Delta_{re}(\D):=W(B)B, \quad B:=\{[E_1],\dots, [E_\mu]\},
\end{equation}
where $W(B)$ is a subgroup of ${\rm Aut}(K_0(\D),I_\D)$ generated by reflections
\begin{equation}
r_{[E_i]}(\lambda):=\lambda-I_\D(\lambda,[E_i])[E_i],\quad \lambda\in K_0(\D),\quad i=1,\dots, \mu, 
\end{equation}
\item 
the automorphism $c_{\D}$ on $K_0(\D)$ induced by the Coxeter functor
$\C_{\D}:=\S_\D[-1]$ on $\D$ where $\S_\D$ is the Serre functor on $\D$,
\end{itemize}
forms a simply-laced generalized root system $R_\D$, which does not depend on the choice of the full exceptional collection $\E$.
\end{prop}
\begin{pf}
Since $\E=(E_1,\dots, E_\mu)$ is a full strongly exceptional collection in $\D$, it follows that $K_0(\D)=\displaystyle\bigoplus_{i=1}^\mu\ZZ[E_i]$ on which the Euler form $\chi_\D$ is well-defined and there exists 
the Serre functor $S_\D$ on $\D$. 
\begin{lem}
We have
\begin{equation}
c_\D=r_{[E_1]}\cdots r_{[E_\mu]}.
\end{equation}
\end{lem}
\begin{pf}
The statement follows from a relation between the Serre functor $\S_\D$ on $\D$ and the helix generated by 
the full exceptional collection $\E$.
See p.~223 in \cite{BP}, for example.
\qed
\end{pf}
We only have to show that $\Delta_{re}(\D)$ satisfies the desired properties and it does not depend on the particular choice of $\E$ since the objects $K_0(\D)$, $I_\D$ and $c_\D$ are invariants of the triangulated category $\D$.
It is obvious from the definition of exceptional object that $I_\D(\alpha,\alpha)=2$ for all $\alpha\in\Delta_{re}(\D)$.
Set
\begin{equation}
W(\D):=\{r_{\alpha}\in {\rm Aut}(K_0(\D),I_\D)~\vert~\alpha\in\Delta_{re}(\D)\}.
\end{equation}
\begin{lem}\label{lem:29}
For any $\alpha\in \Delta_{re}(\D)$, we have
\begin{equation}
r_{[E_i]}r_{\alpha}=r_{r_{[E_i]}(\alpha)}r_{[E_i]}.
\end{equation}
\end{lem}
\begin{pf}
A direct calculation yield the statement.
\qed
\end{pf}
Note that Lemma~\ref{lem:29} implies that $W(\D)=W(B)$.
\begin{lem}\label{lem:212}
For an exceptional object $E'\in \D$, the class $[E']\in K_0(\D)$ belongs to $\Delta_{re}(\D)$.
\end{lem}
\begin{pf}
Choose a full exceptional collection $\E'$ in $\D$ such that $E'\in\E'$. 
There exists an element $g\in G_\mu$ such that $g\E'\cong \E$.
Since we have $-[X]=r_{[X]}([X])$ and $[\R_X Y]=-r_{[Y]}[X]=r_{[Y]}r_{[X]}([X])$ in $K_0(\D)$ for any exceptional pair $(X,Y)$, 
it turns out from Lemme~\ref{lem:29} that $[E']\in\Delta_{re}(\D)$.
\qed
\end{pf}
Set $B':=\{[E'_1],\dots. [E'_\mu]\}$ for any full exceptional collection $\E'=(E'_1,\dots, E'_\mu)$ in $\D$.
Lemma~\ref{lem:212} implies that $W(\D)B'\subset W(\D)W(\D)B\subset W(\D)B$ and hence $W(\D)B'=W(\D)B$. 
Therefore the set $\Delta_{re}(\D)$ does not depend on the particular choice of the full exceptional collection $\E$.
Thus we have completed the proof of the proposition.
\qed
\end{pf}
\begin{rem}
We assumed in Proposition~\ref{prop:1} the existence of a full {\it strongly} exceptional collection $\E$ in $\D$ 
in order to ensure that $\D$ has a unique enhancement in a suitable sense. 
We refer \cite{kajiura} and \cite{lo} for some results on the uniqueness of enhancements for triangulated categories 
and do not discuss this matter more in detail.
\end{rem}
\begin{defn}
The generalized root system $R_\D$ in Proposition~\ref{prop:1} is called the {\it simply-laced generalized root system associated to $\D$}.
\end{defn}
From various points of view, which we do not discuss in detail in this paper, it is natural to expect the assumptions of Proposition~\ref{prop:1}. Indeed, they are proven for derived categories of hereditary Artin algebras by Crawley-Boevey \cite {CB1} and Ringel \cite{ringel} and for derived categories of coherent sheaves on an orbifold projective line $\PP^1_{A,\Lambda}$ (we shall recall the definition later) by Meltzer \cite{meltzer}.
The transitivity of the action of $G_\mu$ is conjectured by Bondal--Polishchuk (Conjecture~2.2 in \cite{BP}), 
and is known for the derived categories of coherent sheaves on $\PP^2$ and $\PP^1\times \PP^1$ by Rudakov \cite{rudakov},  
by arbitrary del Pezzo surfaces by Kuleshov and Orlov \cite{kuleshov-orlov}, for example.
\begin{rem}
One can also consider the subset $\Delta_{re}^{s}(\D)$ of $K_0(\D)$ defined by 
\begin{equation}
\Delta_{re}^{s}(\D):=\{[E]\in K_0(\D)~\vert~E\text{ is an exceptional object in }\D \},
\end{equation}
which is known as the set of {\it Schur roots}. 
Under the assumptions of Proposition~\ref{prop:1}, we always have $\Delta_{re}^{s}(\D)\subset \Delta_{re}(\D)$, however, 
$\Delta_{re}^{s}(\D)\ne \Delta_{re}(\D)$ in general.
A criteria to have $\Delta_{re}^{s}(\D)$ in terms of the Weyl group $W(\D)$ is recently given by Hubery--Krause \cite{HK} 
for derived categories of hereditary Artin algebras.
\end{rem}
\subsection{Generalized root systems associated to star quivers}
We recall the definition of quivers and their path algebras.
\begin{defn}
A {\it quiver} $Q$ is a quadruple $(Q_0,Q_1;s,t)$ where $Q_0$ is 
a set called the set of {\it vertices}, $Q_1$ is 
a set called the set of {\it arrows} and $s,t$ are maps from $Q_1$ to $Q_0$ 
which associate the {\it source} vertex and the {\it target} vertex for each arrow. 
An arrow $f$ with the source $s(f)$ and the target $t(f)$ is often written as $s(f)\stackrel{f}{\longrightarrow }t(f)$.
\end{defn}
\begin{defn}
Let $Q=(Q_0,Q_1;s,t)$ be a quiver.
\begin{enumerate}
\item
A {\it path of length $0$} is a symbol $(v|v)$ defined for each vertex $v\in Q_0$. 
\item
A {\it path of length $l\ge 1$} from the vertex $v$ to the vertex $v'$ in a quiver $Q$ 
is a symbol $(v|f_1\cdots f_l|v')$ with arrows $f_i$, $i=1,\dots, l$ such that $s(f_1)=v$, $t(f_l)=v'$ and 
$s(f_{i+1})=t(f_{i})$, $i=1,\dots, l-1$.
\item
For a path $p=(v|f_1\cdots f_l|v')$, set $s(p):=v$ and $t(p):=v'$. 
\item
An ordered pair of paths $(p_1,p_2)$ is {\it composable} if $t(p_1)=s(p_2)$.
\item
The {\it composition} of composable paths $\left((v_1|f_1\cdots f_l|v_1'), (v_2|g_1\cdots g_m|v_2')\right)$ 
is a path $(v_1|f_1\cdots f_l g_1\cdots g_m|v_2')$.
\end{enumerate}
\end{defn}
\begin{defn}
Let $Q$ be a quiver.
\begin{enumerate}
\item
The {\it path algebra} $k  Q$ of a quiver $Q$ 
is defined as the $k$-module generated by all paths in $Q$ together with 
the associative product structure defined by the composition of paths, 
where the product of two non-composable paths is set to be zero. 
\item
A {\it bound quiver} is a pair $(Q,\I)$ where $Q$ is a quiver and $\I$ is an ideal of $k  Q$.  
\item
A {\it bound quiver algebra} $k (Q,\I)$ of a bound quiver $(Q,\I)$ is defined as the algebra $k  Q/\I$.
\end{enumerate}
\end{defn}
We recall a special class of quivers called star quivers, which are of our interest.
\begin{defn}\label{defn:2.19}
Let $r\ge 3$ be a positive integer and $A = (a_{1},\ldots,a_{r})$ a tuple of positive integers greater than one.
Define a quiver $\TT_{A}=\left(\TT_{A,0},\TT_{A,1};s,t\right)$ as follows$:$
\begin{subequations} 
\begin{itemize}
\item
The set $\TT_{A,0}$ of vertices is  
\begin{equation}
\TT_{A,0} := \{ 1\}\coprod\left(\coprod_{i=1}^r\coprod_{j=1}^{a_i-1}\{(i,j)\}\right).
\end{equation} 
\item
The set $\TT_{A,1}$ of arrows is 
\begin{equation} 
\TT_{A,1} := \coprod_{i=1}^r\coprod_{j=1}^{a_i-1}\{f_{i,j}\},
\end{equation} 
whose source $s(f)$ and target $t(f)$ of each arrow $f$ is given as follows$;$
\begin{equation}
s(f_{i,1})=1,\quad t(f_{i,1})=(i,1),\quad i=1,\ldots, r,
\end{equation}
\begin{equation}
s(f_{i,j})=(i,j-1),\quad t(f_{i,j})=(i,j),\quad i=1,\dots, r,\ j=1,\ldots,a_{i}-1.
\end{equation}
\end{itemize}
\end{subequations}
The quiver $\TT_{A}$ is called the {\it star quiver of type $A$}. 
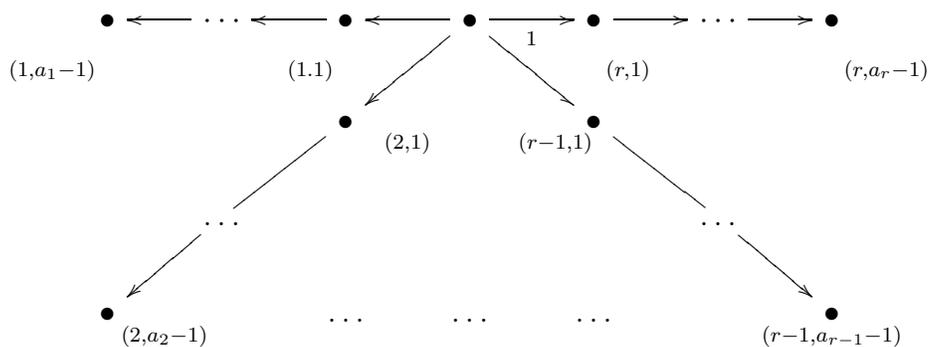
\[
\xymatrix{ 
{\bullet} \ar@{<-}[r] \ar@{}_{(1,a_{1}-1)}[d]  & {\cdots} \ar@{<-}[r]  & {\bullet} \ar@{<-}[r]   \ar@{}_{(1.1)}[d] & {\bullet} \ar@{->}[dl] \ar@
{->}[r] \ar@{}_{1}[r] \ar@_{->}[dr]& {\bullet} \ar@{->}[r]  \ar@{}^{(r,1)}[d]  & {\cdots} \ar@{->}[r]  &{\bullet} \ar@{}^{(r,a_{r}-1)}[d]  \\
& &  {\bullet} \ar@{-}[dl] \ar@{}_{(2,1)}[r]  & &  {\bullet} \ar@{-}[dr] \ar@{}_{(r-1,1)} & &\\
 & {\dots} \ar@{->}[dl] & & & &  {\dots} \ar@{->}[dr] &\\
{\bullet}  \ar@{}_{(2,a_{2}-1)}[r] & & \dots & \dots& \dots & & {\bullet}  \ar@{}_{(r-1,a_{r-1}-1)}
  }
\]
\end{defn}

\begin{defn}\label{defn of TT_A}
Let $\TT_{A}$ be a star quiver of type $A$.
\begin{enumerate}
\item
Denote by $R_A$ the generalized root system associated to $\D^b(k\TT_A)$.
\item 
Let $\alpha_v$ be the equivalence class in $K_0(R_A)=K_0(\D^b(k\TT_A))$ of the simple $k\TT_A$-module corresponding to the vertex $v\in \TT_{A,0}$.
Set
\begin{equation}
B_{\TT_A}:=\{\alpha_v\}_{v\in\TT_{A,0}},
\end{equation}  
which is a root basis of $R_A$. 
\item
Denote by $T_A$ the Coxeter--Dynkin diagram for $\Gamma_{B_{\TT_{A}}}$, which is given by 
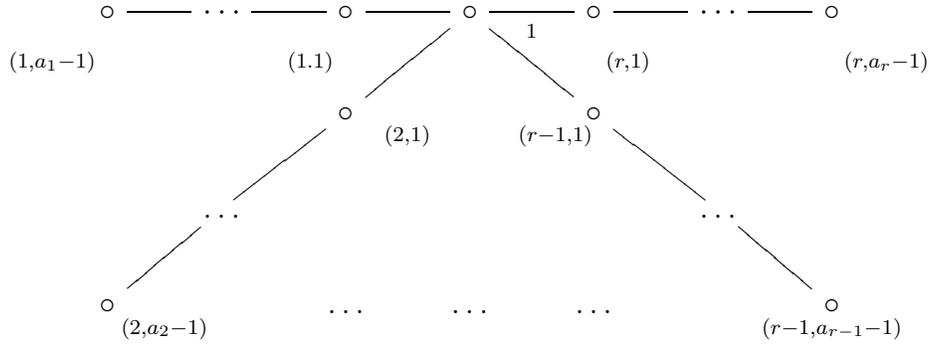
\[
\xymatrix{ 
{\circ} \ar@{-}[r] \ar@{}_{(1,a_{1}-1)}[d]  & {\cdots} \ar@{-}[r]  & {\circ} \ar@{-}[r]    \ar@{}_{(1.1)}[d] & {\circ} \ar@{-}[dl] \ar@
{-}[r] \ar@{}_{1}[r] \ar@_{-}[dr]& {\circ} \ar@{-}[r]  \ar@{}^{(r,1)}[d]  & {\cdots} \ar@{-}[r]  &{\circ} \ar@{}^{(r,a_{r}-1)}[d]  \\
& &  {\circ} \ar@{-}[dl] \ar@{}_{(2,1)}[r]  & &  {\circ} \ar@{-}[dr] \ar@{}_{(r-1,1)} & &\\
 & {\dots} \ar@{-}[dl] & & & &  {\dots} \ar@{-}[dr] &\\
{\circ}  \ar@{}_{(2,a_{2}-1)}[r] & & \dots & \dots& \dots & & {\circ}  \ar@{}_{(r-1,a_{r-1}-1)}
  }
\]
We often write $v\in {T}_A$ instead of $v\in{\TT}_{A,0}$.
\item
For each $v\in T_{A}$, define the {\it simple reflection} $r_v$ on $K_0(R_{A})_\QQ$ by
\begin{equation}
r_v(\lambda):=\lambda-I_{R_A}(\lambda ,\alpha_v)\alpha_v,\quad \lambda\in K_0(R_{A})_\QQ.
\end{equation}
Since $B_{\TT_A}$ is a root basis of $R_A$, the Weyl group $W(R_A)$ of $R_A$ is generated by simple reflections;
\begin{equation}
W(R_A)=\langle r_v~\vert~v\in T_{A}\rangle.
\end{equation}
\end{enumerate}
\end{defn}
Note that the Cartan matrix $\left(I_{R_A}(\alpha_v,\alpha_{v'})\right)$ is a generalized Cartan matrix in the sense of \cite{Kac1}.
Therefore one can naturally associate to $R_A$ a Kac--Moody Lie algebra $\g(R_A)$.
\subsection{Octopus}
We introduce a bound quiver, a ``one point extension" of the star quiver. 
\begin{defn}\label{defn of widetilde{TT}_A}
Let $r\ge 3$ be a positive integer, $A = (a_{1},\ldots,a_{r})$ an  $r$-tuple of positive integers greater than one and 
$\Lambda=(\lambda_1,\ldots,\lambda_r)$ an $r$-tuple of pairwise distinct elements of $\PP^1(k)$ normalized such that 
$\lambda_1=\infty$, $\lambda_2=0$ and $\lambda_3=1$.
\begin{enumerate}
\item
Define a quiver $\widetilde{\TT}_A=(\widetilde{\TT}_{A,0},\widetilde{\TT}_{A,1},s,t)$ as follows$:$
\begin{subequations} 
\begin{itemize}
\item
The set $\widetilde{\TT}_{A,0}$ of vertices is given by 
\begin{equation}
\widetilde{\TT}_{A,0} := {\TT}_{A,0}\coprod \{1^*\}=
\{ 1\}\coprod\left(\coprod_{i=1}^r\coprod_{j=1}^{a_i-1}\{(i,j)\}\right)\coprod\{1^*\}.
\end{equation} 
\item
The set $\widetilde{\TT}_{A,1}$ of arrows is given by 
\begin{equation} 
\widetilde{\TT}_{A,1} := \TT_{A,1}\coprod\left(\coprod_{i=1}^r\{f_{i,1^*}\}\right)=
\left(\coprod_{i=1}^r\coprod_{j=1}^{a_i-1}\{f_{i,j}\}\right)\coprod \left(\coprod_{i=1}^r\{f_{i,1^*}\}\right),
\end{equation} 
whose source $s(f)$ and target $t(f)$ of each arrow $f$ is given as follows:
\begin{equation}
s(f_{i,1})=1,\quad t(f_{i,1})=(i,1),\quad i=1,\ldots, r,
\end{equation}
\begin{equation}
s(f_{i,j})=(i,j-1),\quad t(f_{i,j})=(i,j),\quad i=1,\dots, r,\ j=2,\ldots,a_{i}-1,
\end{equation}
\begin{equation}
s(f_{i,1^*})=(i,1),\quad t(f_{i,1^*})=1^*,\quad i=1,\dots, r.
\end{equation}
\end{itemize}
\item
Define an ideal $\I_\Lambda$ of the path algebra $k\TT_A$ by
\begin{equation} 
\I_\Lambda:=\left<
\sum_{i=1}^r\lambda_i^{(1)}f_{i,1}f_{i,1^*}, 
\sum_{i=1}^r\lambda_i^{(2)}f_{i,1}f_{i,1^*}
\right>,
\end{equation} 
\end{subequations}
where $(\lambda_1^{(1)},\lambda_1^{(2)})=(1,0)$ and $(\lambda_i^{(1)},\lambda_i^{(2)})=(\lambda_i,1)$ for 
$i=2,\dots ,r$.
\end{enumerate}
We denote by $\widetilde{\TT}_{A,\Lambda}$ the bound quiver $(\widetilde{\TT}_A,\I_\Lambda)$ for simplicity.
The bound quiver algebra $k\widetilde{\TT}_{A,\Lambda}$ is called the {\it octopus} of type $(A,\Lambda)$.
\[
\xymatrix{ 
 & & & {\bullet}  \ar@{<-}[dr]  \ar@{<-}[ldd] \ar@2{.}[d] \ar@{<-}[rdd] \ar@{}^{1^{*}}[r] & & &  \\
{\bullet} \ar@{<-}[r] \ar@{}_{(1,a_{1}-1)}[d]  & {\cdots} \ar@{<-}[r]  & {\bullet} \ar@{<-}[r] \ar@{->}[ur]   \ar@{}_{(1.1)}[d] & {\bullet} \ar@{->}[dl] \ar@
{->}[r] \ar@{}_{1}[r] \ar@_{->}[dr]& {\bullet} \ar@{->}[r]  \ar@{}^{(r,1)}[d]  & {\cdots} \ar@{->}[r]  &{\bullet} \ar@{}^{(r,a_{r}-1)}[d]  \\
& &  {\bullet} \ar@{->}[dl] \ar@{}_{(2,1)}[r]  & &  {\bullet} \ar@{->}[dr] \ar@{}_{(r-1,1)} & &\\
 & {\dots} \ar@{->}[dl] & & & &  {\dots} \ar@{->}[dr] &\\
{\bullet}  \ar@{}_{(2,a_{2}-1)}[r] & & \dots & \dots& \dots & & {\bullet}  \ar@{}_{(r-1,a_{r-1}-1)}
  }
\]
\end{defn}
\begin{rem}
In \cite{CB2}, Clawley-Boevey defines a bound quiver algebra associated to $(A,\Lambda)$, which is called the {\it squid}.
A squid and an octopus are different but very similar, more precisely, these algebras are not isomorphic but derived equivalent.
\end{rem}
\subsection{Algebro-geometric aspect of octopuses}
We associate to a pair $(A,\Lambda)$ an algebro-geometric object following Geigle--Lenzing (cf. Section~1.1 in \cite{gl:1}).
\begin{defn}
Let $r\ge 3$ be a positive integer, $A = (a_{1},\ldots,a_{r})$ an  $r$-tuple of positive integers greater than one and 
$\Lambda=(\lambda_1,\ldots,\lambda_r)$ an $r$-tuple of pairwise distinct elements of $\PP^1(k)$ normalized such that 
$\lambda_1=\infty$, $\lambda_2=0$ and $\lambda_3=1$.
\begin{enumerate}
\item
Define a ring $S_{A,\Lambda}$ by 
\begin{equation} 
S_{A,\Lambda}:=k[X_1,\dots,X_r]\left/(X_i^{a_i}-X_2^{a_2}+\lambda_i X_1^{a_1};i=3,\dots, r)\right..
\end{equation}
\item
Denote by $L_A$ an abelian group generated by $r$-letters $\vec{X_i}$, $i=1,\dots ,r$ defined as the quotient 
\begin{equation}
L_A:=\bigoplus_{i=1}^r\ZZ\vec{X}_i\left/(a_i\vec{X}_i-a_j\vec{X}_j;1\le i<j\le r)\right. .
\end{equation}
\end{enumerate}
\end{defn}
Note that $S_{A,\Lambda}$ is naturally graded with respect to $L_A$.
Denote by ${\rm gr}^{L_A}\text{-}S_{A,\Lambda}$ the category of finitely generated 
$L_A$-graded $S_{A,\Lambda}$-modules and 
by ${\rm tor}^{L_A}\text{-}S_{A,\Lambda}$ the full subcategory of ${\rm gr}^{L_A}\text{-}S_{A,\Lambda}$ 
consisting of modules of finite length.
\begin{defn}\label{def:2.22}
Define a stack $\PP^1_{A,\Lambda}$ by
\begin{equation}
\PP^1_{A,\Lambda}:=\left[\left({\rm Spec}(S_{A,\Lambda})\backslash\{0\}\right)/{\rm Spec}({k L_A})\right],
\end{equation}
which is called the {\it orbifold projective line} of type $(A,\Lambda)$. 
Denote by ${\rm coh}(\PP^1_{A,\Lambda})$ the category of coherent sheaves on $\PP^1_{A,\Lambda}$ and 
by $\D^b{\rm coh}(\PP^1_{A,\Lambda})$ its bounded derived category.
\end{defn}
Properties of categories ${\rm coh}(\PP^1_{A,\Lambda})$ and $\D^b{\rm coh}(\PP^1_{A,\Lambda})$ are extensively studied by 
Geigle--Lenzing \cite{gl:1}. Among them, the following is of our interest in this paper.
\begin{prop}\label{prop:gl}
There exists an equivalence of triangulated categories
\begin{equation}
\D^b{\rm coh}(\PP^1_{A,\Lambda})\simeq \D^b(k\widetilde{\TT}_{A,\Lambda}).
\end{equation}
\end{prop}
\begin{pf}
It follows from Proposition~4.1 in \cite{gl:1}.
\qed
\end{pf}

\subsection{Generalized root systems associated to octopuses}
Since the assumptions of Proposition~\ref{prop:1} are proven for $\D^b{\rm coh}(\PP^1_{A,\Lambda})$ by Meltzer \cite{meltzer},
we obtain a generalized root system.
\begin{defn}\label{def:226}
Let $k\widetilde{\TT}_{A,\Lambda}$ be an octopus of type $(A,\Lambda)$.
\begin{enumerate}
\item
Denote by $\widetilde{R}_A$ the generalized root system associated to
$\D^b(k\widetilde{\TT}_{A,\Lambda})$.
\item
For any $v\in \widetilde{\TT}_{A,0}$, denote by $P_v$ the corresponding
indecomposable projective $k\widetilde{\TT}_{A,\Lambda}$-module, which
satisfies $k\widetilde{\TT}_{A,\Lambda}=\displaystyle
\bigoplus_{v\in\widetilde{T}_A}P_v$ as a
$k\widetilde{\TT}_{A,\Lambda}$-module.
Note that the collection $(P_v)_{v\in\widetilde{T}_A}$ forms a full
strongly exceptional collection in $\D^b(k\widetilde{\TT}_{A,\Lambda})$.
\item
For any $v\in \widetilde{\TT}_{A,0}$, denote by $S_v$ the corresponding
simple $k\widetilde{\TT}_{A,\Lambda}$-module.
Note that the collection $(S_v)_{v\in\widetilde{T}_A}$ forms a full
exceptional collection in $\D^b(k\widetilde{\TT}_{A,\Lambda})$
such that
\begin{equation}\label{halfdual}
\chi_{\D^b(k\widetilde{\TT}_{A,\Lambda})}([S_v],[P_{v'}])=\delta_{vv'},\quad
v,v'\in \widetilde{\TT}_{A,0},
\end{equation}
where $\delta_{vv'}$ denotes the Kronecker's delta.
\item
For any simple $k\widetilde{\TT}_{A,\Lambda}$-module $S_v$, $v\in
\widetilde{\TT}_{A,0}$, denote by
$\widetilde{\alpha}_v$ the equivalence class $[S_v]\in
K_0(\widetilde{R}_A)=K_0(\D^b(k\widetilde{\TT}_{A,\Lambda}))$.
Set
\begin{equation}
B_{\widetilde{\TT}_{A,\Lambda}}:=\{\widetilde{\alpha}_v\}_{v\in\widetilde{\TT}_{A,0}},
\end{equation}
which is a root basis of $\widetilde{R}_A$.
\item
Denote by $\widetilde{T}_A$ the Coxeter--Dynkin diagram
$\Gamma_{B_{\widetilde{\TT}_{A,\Lambda}}}$, which turns out to be
the following diagram by using the property~\eqref{halfdual}:
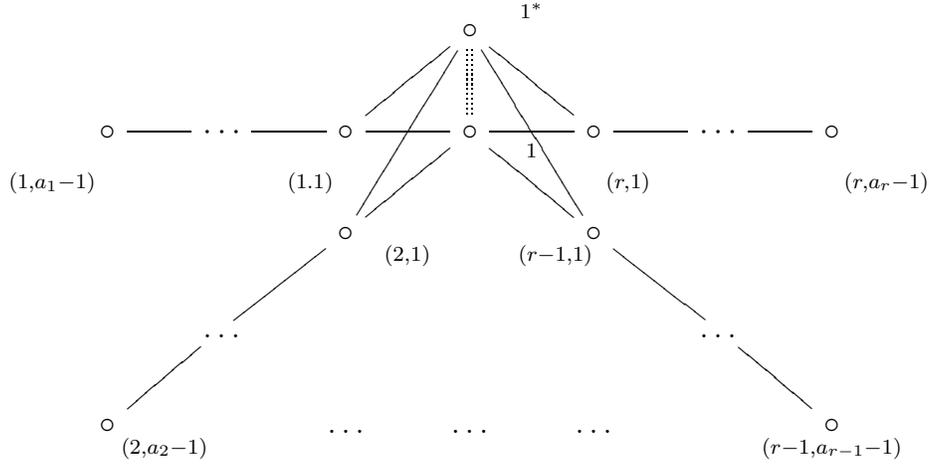
\[
\xymatrix{
 & & & {\circ}  \ar@{-}[dr]  \ar@{-}[ldd] \ar@2{.}[d] \ar@{-}[rdd]
\ar@{}^{1^{*}}[r] & & &  \\
{\circ} \ar@{-}[r] \ar@{}_{(1,a_{1}-1)}[d]  & {\cdots} \ar@{-}[r]  &
{\circ} \ar@{-}[r] \ar@{-}[ur]   \ar@{}_{(1.1)}[d] & {\circ} \ar@{-}[dl]
\ar@
{-}[r] \ar@{}_{1}[r] \ar@_{-}[dr]& {\circ} \ar@{-}[r]  \ar@{}^{(r,1)}[d] 
& {\cdots} \ar@{-}[r]  &{\circ} \ar@{}^{(r,a_{r}-1)}[d]  \\
& &  {\circ} \ar@{-}[dl] \ar@{}_{(2,1)}[r]  & &  {\circ} \ar@{-}[dr]
\ar@{}_{(r-1,1)} & &\\
 & {\dots} \ar@{-}[dl] & & & &  {\dots} \ar@{-}[dr] &\\
{\circ}  \ar@{}_{(2,a_{2}-1)}[r] & & \dots & \dots& \dots & & {\circ} 
\ar@{}_{(r-1,a_{r-1}-1)}
  }
\]
We often write $v\in\widetilde{T}_A$ instead of $v\in\widetilde{\TT}_{A,0}$.
\item
For each $v\in \widetilde{T}_{A}$, define the {\it simple reflection}
$\widetilde{r}_v$ on $K_0(\widetilde{R}_A)$ by
\begin{equation}
\widetilde{r}_v(\widetilde{\lambda}):=\widetilde{\lambda}-I_{\widetilde{R}_A}(\widetilde{\lambda}
,\widetilde{\alpha}_v)\widetilde{\alpha}_v,\quad \widetilde{\lambda}\in
K_0(\widetilde{R}_A).
\end{equation}
Since $B_{\widetilde{\TT}_A}$ is a root basis of $\widetilde{R}_A$, the
Weyl group $W(\widetilde{R}_A)$ of $\widetilde{R}_A$ is generated by
simple reflections;
\begin{equation}
W(\widetilde{R}_A)=\langle \widetilde{r}_v~\vert~v\in
\widetilde{T}_{A}\rangle.
\end{equation}
\end{enumerate}
\end{defn}
\subsection{A relation between octopuses and star quivers}\label{subsec 2.7}
Set $\delta:=\widetilde{\alpha}_{1^*}-\widetilde{\alpha}_{1}$.
It is easy to see that $\delta$ belongs to the radical of the Cartan form $I_{\widetilde{R}_A}$ on $K_0(\widetilde{R}_A)$, 
therefore the natural projection map 
\begin{equation}\label{14}
K_0(\widetilde{R}_A)\longrightarrow K_0(\widetilde{R}_A)/\ZZ\delta \cong K_0(R_A)
\end{equation}
induces the surjective group homomorphism 
\begin{equation}
p:W(\widetilde{R}_A)\twoheadrightarrow  W(R_A).
\end{equation}
Indeed, we have
\begin{subequations}
\begin{equation}
p(\widetilde{r}_{1})=p(\widetilde{r}_{{1^*}})=r_{1},
\end{equation}
\begin{equation}
p(\widetilde{r}_{v})=r_{v},\quad v\in \TT_{A,0}.
\end{equation}
\end{subequations}
Moreover, the correspondence $\alpha_v\mapsto \widetilde{\alpha}_v$ for $v\in T_{A}$ gives the splitting 
of the surjective map \eqref{14} and induces the isomorphism of $\ZZ$-modules
\begin{equation}
K_0(\widetilde{R}_A)\cong K_0(R_A)\oplus \ZZ\delta,
\end{equation}
which is compatible with the Cartan forms $I_{\widetilde{R}_A}$ and $I_{R_A}$. 
Hence we obtain the group homomorphism
\begin{equation}
i:W(R_A)\longrightarrow W(\widetilde{R}_A),\quad r_{v}\mapsto \widetilde{r}_{v}
\end{equation}
such that $p\circ i={\rm id}_{W(R_A)}$.
\section{Presentations of Weyl groups}\label{Presentations of Weyl groups}
In this section, we describe the Weyl group $W(\widetilde{R}_A)$ as the ``affinization" of the Weyl group $W(R_A)$. 
\begin{defn}
For each vertex $v\in T_A$, define an element $\widetilde{\tau}_{v}\in W(\widetilde{R}_A)$ by induction as follows$:$
\begin{itemize}
\item
For the vertex $1$, set
\begin{subequations}
\begin{equation}
\widetilde{\tau}_{1}:=\widetilde{r}_{1}\widetilde{r}_{{1^*}}.
\end{equation}
\item
Set
\begin{equation}
\widetilde{\tau}_{(i,1)}:=\widetilde{r}_{{(i,1)}}\widetilde{\tau}_{1}\widetilde{r}_{{(i,1)}}\widetilde{\tau}_{1}^{-1},
\quad i=1,\dots, r,
\end{equation}
\begin{equation}
\widetilde{\tau}_{(i,j)}:=\widetilde{r}_{{(i,j)}}\widetilde{\tau}_{(i,j-1)}\widetilde{r}_{{(i,j)}}\widetilde{\tau}_{(i,j-1)}^{-1},
\quad i=1,\dots, r,\ j=2,\dots, a_i-1.
\end{equation}
\end{subequations}
\end{itemize}
\end{defn}
Denote by $N$ the smallest normal subgroup of $W(\widetilde{R}_A)$ containing $\widetilde{\tau}_{1}$.
\begin{lem}\label{normal N gamma}
For all $v\in T_A$, the element $\widetilde{\tau}_{v}$ belongs to $N$.
\end{lem}
\begin{pf}
By the fact that $\widetilde{r}_{{v}}^2=1$ and the definition of $N$, 
we have $\widetilde{r}_{{(i,1)}}\widetilde{\tau}_{1}\widetilde{r}_{{(i,1)}} \in N$ and hence $\widetilde{\tau}_{(i,1)} \in N$
for $i=1,\dots, r$. We shall show this lemma by the induction on $j$. We assume that $\widetilde{\tau}_{(i,j)} \in N$. Under this assumption,
one has $\widetilde{r}_{{(i,j+1)}}\widetilde{\tau}_{(i,j)}\widetilde{r}_{{(i,j+1)}} \in N$ and hence $\widetilde{\tau}_{(i,j+1)} \in N$.
Therefore this lemma holds.
\qed
\end{pf}

\begin{prop}\label{translation}
For all $v\in T_A$, we have
\begin{equation}\label{34}
\widetilde{\tau}_{v}(\widetilde{\lambda})=\widetilde{\lambda}-I_{\widetilde{R}_A}(\widetilde{\lambda}, \widetilde{\alpha}_v)\delta,
\quad \widetilde{\lambda}\in K_0(\widetilde{R}_A).
\end{equation}
In particular, there is a natural surjective group homomorphism 
\begin{equation}
\varphi: K_0(R_A)\twoheadrightarrow  N,\quad \sum_{v\in T_A}m_v \alpha_v
\mapsto \prod_{v\in T_A}\widetilde{\tau}_{v}^{m_v},
\end{equation}
which induces an isomorphism
\begin{equation}
K_0(R_A)/{\rm rad}(I_{R_A})\cong N.
\end{equation}
\end{prop}
\begin{pf}
We have
\begin{eqnarray*}
\widetilde{\tau}_{1}(\widetilde{\lambda})
&=& \widetilde{r_{1}}\widetilde{r}_{1^{*}}(\widetilde{\lambda}) \\
&=& \widetilde{r_{1}} (\widetilde{\lambda} - I_{\widetilde{R}_A}(\widetilde{\lambda} , \widetilde{\alpha}_{1^{*}}) \widetilde{\alpha}_{1^{*}}) \\
&=& \widetilde{\lambda} - I_{\widetilde{R}_A}(\widetilde{\lambda} , \widetilde{\alpha}_{1}) (\widetilde{\alpha}_{1^{*}} - \widetilde{\alpha}_{1}) 
=\widetilde{\lambda} - I_{\widetilde{R}_A}(\widetilde{\lambda} , \widetilde{\alpha}_{1}) \delta.
\end{eqnarray*}
We have a following formula in the same way:
\begin{equation*}
\widetilde{\tau}_{(i,1)}(\widetilde{\lambda}) = \widetilde{\lambda} - I(\widetilde{\lambda} , \widetilde{\alpha}_{(i,1)}) \delta, 
\quad i=1,\dots, r.
\end{equation*}
We show the equation \eqref{34} by the induction for $i$.
Assume the following expression:
\begin{equation*}
\widetilde{\tau}_{(i,k)} (\widetilde{\lambda})=\widetilde{\lambda} - I_{\widetilde{R}_A}(\widetilde{\lambda} , \widetilde{\alpha}_{(i,k)}) \delta,\quad k=2,\dots, a_{i}-2.
\end{equation*}
Then we have
\begin{eqnarray*}
\widetilde{\tau}_{(i,k+1)} 
&=& \widetilde{r}_{(i,k+1)} \widetilde{\tau}_{(i,k)} \widetilde{r}_{(i,k+1)} \widetilde{\tau}_{(i,k)}^{-1} (\widetilde{\lambda}) \\
&=& \widetilde{r}_{(i,k+1)} \widetilde{\tau}_{(i,k)} \widetilde{r}_{(i,k+1)} (\widetilde{\lambda} + I_{\widetilde{R}_A}(\widetilde{\lambda} , \widetilde{\alpha}_{(i,k)}) \widetilde{\alpha}_{(i,k)}) \\
&=&  \widetilde{r}_{(i,k+1)} \widetilde{\tau}_{(i,k)} ( \widetilde{\lambda} - I_{\widetilde{R}_A}(\widetilde{\lambda} , \widetilde{\alpha}_{(i,k+1)})\widetilde{\alpha}_{(i,k+1)} + I_{\widetilde{R}_A}(\widetilde{\lambda} , \widetilde{\alpha}_{(i,k)}) \delta )  \\
&=& \widetilde{r}_{(i,k+1)} (\widetilde{\lambda} - I_{\widetilde{R}_A}(\widetilde{\lambda} , \widetilde{\alpha}_{(i,k+1)}) \widetilde{\alpha}_{(i,k+1)} - I_{\widetilde{R}_A}(\widetilde{\lambda} , \widetilde{\alpha}_{(i,k+1)}) \delta) \\
&=& \widetilde{\lambda} - I_{\widetilde{R}_A}(\widetilde{\lambda} , \widetilde{\alpha}_{(i,k+1)}) \delta.
\end{eqnarray*}
Hence we obtain the equation \eqref{34}.
We calculate the kernel of the natural surjective group homomorphism $\varphi$.
We have
\begin{equation*}
\prod_{v\in T_{A}} \widetilde{\tau}_{v}^{m_{v}} (\widetilde{\lambda}) =\widetilde{\lambda} - \sum_{v\in T_{A}} m_{v} I_{\widetilde{R}_A}(\widetilde{\lambda} , \widetilde{\alpha}_{v}) \delta.
\end{equation*}
An element $\displaystyle\sum_{v\in T_{A}} m_{v} {\alpha}_{v} \in K_0(R_A)$ is in ${\rm Ker}(\varphi)$ if and only if $\displaystyle I_{\widetilde{R}_A}\left(\widetilde{\lambda} , \sum_{v\in T_{A}} m_{v} \widetilde{\alpha}_{v}\right) \delta = 0$ for all $\widetilde{\lambda}$, which means $\displaystyle\sum_{v\in T_{A}} m_{v} \widetilde{\alpha}_{v}\in {\rm rad}(I_{R_A})$.
\qed
\end{pf}
Note that ${\rm rad}(I_{R_A})$ is zero if $\chi_A\ne 0$ and is of rank one if $\chi_A=0$.
\begin{prop}\label{prop:34}
For $v,$ $v'\in {T}_{A}$, we have
\begin{subequations}
\begin{equation}\label{T-adj0}
\widetilde{r}_{v}\widetilde{\tau}_{v}\widetilde{r}_{v} = \widetilde{\tau}_{v}^{-1}, 
\end{equation}
\begin{equation}\label{T-adj1}
\widetilde{r}_{v}\widetilde{\tau}_{v'}\widetilde{r}_{v} = \widetilde{\tau}_{v'}, \quad \text{if} \quad I_{\widetilde{R}_A}(\widetilde{\alpha}_{v}, \widetilde{\alpha}_{v'}) = 0,
\end{equation}
\begin{equation}\label{T-adj2}
	\widetilde{r}_{v}\widetilde{\tau}_{v'}\widetilde{r}_{v} = \widetilde{\tau}_{v} \widetilde{\tau}_{v'}, \quad \text{if} \quad I_{\widetilde{R}_A}(\widetilde{\alpha}_{v}, \widetilde{\alpha}_{v'})=-1.
\end{equation}
\end{subequations}
\end{prop}
\begin{pf}
First we shall show the equation \eqref{T-adj0}. We have
\begin{eqnarray*}
\widetilde{r}_{v} \widetilde{\tau}_{v} \widetilde{r}_{v} (\widetilde{\lambda}) 
&=& \widetilde{r}_{v} \widetilde{\tau}_{v} (\widetilde{\lambda} - I_{\widetilde{R}_A}(\widetilde{\lambda} , \widetilde{\alpha}_{v})\widetilde{\alpha}_{v})\\
&=& \widetilde{r}_{v} (\widetilde{\lambda} + I_{\widetilde{R}_A}(\widetilde{\lambda} , \widetilde{\alpha}_{v})\delta - I_{\widetilde{R}_A}(\widetilde{\lambda}, \widetilde{\alpha}_{v}) \widetilde{\alpha}_{v} )\\
&=& \widetilde{\lambda} + I_{\widetilde{R}_A}(\widetilde{\lambda} , \widetilde{\alpha}_{v}) \delta= \widetilde{\tau}_{v}^{-1}(\widetilde{\lambda}).
\end{eqnarray*}
Second we shall show the equation \eqref{T-adj1}. We have
\begin{eqnarray*}
\widetilde{r}_{v} \widetilde{\tau}_{v'} \widetilde{r}_{v} (\widetilde{\lambda}) 
&=& \widetilde{r}_{v} \widetilde{\tau}_{v'} (\widetilde{\lambda} - I_{\widetilde{R}_A}(\widetilde{\lambda} , \widetilde{\alpha}_{v})\widetilde{\alpha}_{v})\\
&=& \widetilde{r}_{v} (\widetilde{\lambda} - I_{\widetilde{R}_A}(\widetilde{\lambda} , \widetilde{\alpha}_{v'})\delta - I_{\widetilde{R}_A}(\widetilde{\lambda}, \widetilde{\alpha}_{v}) \widetilde{\alpha}_{v} )\\
&=& \widetilde{\lambda} - I_{\widetilde{R}_A}(\widetilde{\lambda} , \widetilde{\alpha}_{v'}) \delta= \widetilde{\tau}_{v'}(\widetilde{\lambda}).
\end{eqnarray*}
Third we shall show the equation \eqref{T-adj2}. We have
\begin{eqnarray*}
\widetilde{r}_{v} \widetilde{\tau}_{v'} \widetilde{r}_{v} (\widetilde{\lambda}) 
&=& r_{v} \widetilde{\tau}_{v'} (\widetilde{\lambda} - I_{\widetilde{R}_A}(\widetilde{\lambda}, \widetilde{\alpha}_{v})\widetilde{\alpha}_{v})\\
&=& r_{v} (\widetilde{\lambda} - I_{\widetilde{R}_A}(\widetilde{\lambda} , \widetilde{\alpha}_{v'})\delta - I_{\widetilde{R}_A}(\widetilde{\lambda},\widetilde{\alpha}_{v}) \widetilde{\alpha}_{v} - I_{\widetilde{R}_A}(\widetilde{\lambda}, \widetilde{\alpha}_{v})\delta)\\
&=& \widetilde{\lambda} - I_{\widetilde{R}_A}(\widetilde{\lambda}, \widetilde{\alpha}_{v'}) \delta - I_{\widetilde{R}_A}(\widetilde{\lambda}, \widetilde{\alpha}_{v})\delta\\
&=& \widetilde{\tau}_{v+v'}(\widetilde{\lambda})=\widetilde{\tau}_{v}\widetilde{\tau}_{v'}(\widetilde{\lambda}).
\end{eqnarray*}
We have finished the proof of the proposition.
\qed
\end{pf}
Since the Weyl group $W(R_A)$ is a subgroup of ${\rm Aut}(K_0(R_A),I_{R_A})$, we can consider the group $W(R_A)\ltimes K_0(R_A)$, the semi-direct product of $W(R_A)$ and $K_0(R_A)$. 
Note that the equations \eqref{T-adj0}, \eqref{T-adj1} and \eqref{T-adj2} can be thought of
as the adjoint action of $W(R_A)$ on the free generators of $K_0(R_A)$ expressed in multiplicative notation since 
we have $\widetilde{r}_{v}(\widetilde{\alpha}_{v})=-\widetilde{\alpha}_{v}$, $\widetilde{r}_{v}(\widetilde{\alpha}_{v'})=\widetilde{\alpha}_{v'}$ if $I_{\widetilde{R}_A}(\widetilde{\alpha}_{v}, \widetilde{\alpha}_{v'}) = 0$ and 
$\widetilde{r}_{v}(\widetilde{\alpha}_{v'})=\widetilde{\alpha}_{v}+\widetilde{\alpha}_{v'}$ if $I_{\widetilde{R}_A}(\widetilde{\alpha}_{v}, \widetilde{\alpha}_{v'}) = -1$.
Moreover, since the Weyl group $W(R_A)$ respects the radical ${\rm rad}(I_{R_A})$,
we can also consider the group $W(R_A)\ltimes \left(K_0(R_A)/{\rm
rad}(I_{R_A})\right)$, the semi-direct product of $W(R_A)$ and
$K_0(R_A)/{\rm rad}(I_{R_A})$, which is isomorphic to $W(\widetilde{R}_A)$. 
More precisely, we have the following.
\begin{thm}\label{relation W(R_A) and W(widetilde{R}_A)}
There is an exact sequence of groups
\begin{equation}
\{1\}\longrightarrow N\longrightarrow
W(\widetilde{R}_A)\stackrel{p}{\longrightarrow} W(R_A)\longrightarrow
\{1\}.
\end{equation}
In particular, we have an isomorphism
\begin{equation}
W(\widetilde{R}_A)\cong W(R_A)\ltimes \left(K_0(R_A)/{\rm
rad}(I_{R_A})\right).
\end{equation}
\end{thm}
\begin{pf}
It is obvious from Proposition~\ref{translation} that $N$ is a subgroup of the kernel of $p$.
Therefore the natural homomorphism 
\[
i(W(R_A))\ltimes N\longrightarrow W(\widetilde{R}_A),\quad 
(\widetilde{r}_{v}, \widetilde{\tau}_{\alpha_{v'}})\mapsto \widetilde{r}_{v}\widetilde{\tau}_{v'},
\]
is injective. 
On the other hand, the equality
\[
\widetilde{r}_{{1^*}}=\widetilde{r}_{1}\widetilde{r}_{1}\widetilde{r}_{{1^*}}
=\widetilde{r}_{1}\widetilde{\tau}_{1}
\]
implies that this homomorphism is also surjective.
\qed
\end{pf}
Therefore it turns out that $W(\widetilde{R}_A)$ is an affine Weyl group if $\chi_A>0$.
\begin{defn}\label{defn:cuspidal Weyl}
Let the notations be as above.
\begin{enumerate}
\item
If $\chi_A<0$, then the group $W(\widetilde{R}_A)$ is called the {\it cuspidal Weyl group} of type $A$, 
which is isomorphic to $W(R_A)\ltimes K_0(R_A)$ by Theorem~\ref{relation W(R_A) and W(widetilde{R}_A)}.
\item
If $\chi_A=0$, then the group $W(\widetilde{R}_A)$ is called the {\it elliptic Weyl group} of type $A$.
\item
If $\chi_A=0$, then the group $W(R_A)\ltimes K_0(R_A)$ is isomorphic to the non-trivial central extension of $W(\widetilde{R}_A)$ by $\ZZ$, which is called the {\it hyperbolic extension} of the elliptic Weyl group $W(\widetilde{R}_A)$ (cf. Section~1.18 in \cite{saito}).
\end{enumerate}
\end{defn}
\begin{prop}
Under the isomorphism $K_{0}(\widetilde{R}_{A})\cong K_{0}(R_{A})\oplus \ZZ\delta$ in Subsection \ref{subsec 2.7},
we have an isomorphism of sets $\Delta_{re}(\widetilde{R}_A)\cong \Delta_{re}(R_A)+\ZZ\delta$.
\end{prop}
\begin{pf}
In the proof, we identify $K_{0}(\widetilde{R}_{A})$ and $K_{0}(R_{A})\oplus \ZZ\delta$ for simplicity. 
In particular, we have $\widetilde{\alpha}_v=\alpha_v$ for $v\in T_A$ and $\widetilde{\alpha}_{1^{*}}={\alpha}_{1}+\delta$.
\begin{lem}
An element $\alpha \in \Delta_{re}(\widetilde{R}_A)$ can be written in the form $\beta+n\delta$ by 
some $\beta \in \Delta_{re}(R_A)$ and an integer $n\in \ZZ$.
\end{lem}
\begin{pf}
Since any real root is in the orbit of a simple root under the Weyl group action,
a direct calculation yields the following formulas:
\begin{subequations}
\begin{equation}\label{eq:3.8a}
\widetilde{r}_{v}(\beta+n\delta)=\widetilde{r}_{v}(\beta)+n\delta,
\end{equation}
\begin{equation}
\widetilde{r}_{1^{*}}(\beta+n\delta)=\widetilde{r}_{1}(\beta)+I_{\widetilde{R}_A}(\beta, \widetilde{\alpha}_{1})\delta,
\end{equation}
\end{subequations}
where $\beta\in \Delta_{re}(R_{A})$ and $v\in {T}_{A}$. Recall again that $\widetilde{\alpha}_{1^{*}}={\alpha}_{1}+\delta$.
Therefore the iterating simple reflections of a simple root in
$B_{\widetilde{\TT}_{A,\Lambda}}=\{\widetilde{\alpha}_v\}_{v\in T_A}\cup\{\widetilde{\alpha}_{1^{*}}\}$ can be written as the form $\beta+n\delta$ by some $\beta \in \Delta_{re}(R_A)$
and an integer $n\in \ZZ$. 
\qed
\end{pf}
\begin{lem}
For an element $\beta \in \Delta_{re}(R_A)$ and an integer $n\in \ZZ$,
the element $\beta+n\delta \in K_{0}(\widetilde{R}_A)$ is in $\Delta_{re}(\widetilde{R}_A)$.
\end{lem}
\begin{pf}
Recall that $B_{\TT_A}=\{\alpha_v\}_{v\in T_A}$ is a root basis for $R_{A}$.
Since $\Delta_{re}(R_A)=W(R_A)B_{\TT_A}$, we have $\beta+n\delta=r({\alpha}_{v})+n\delta$ for some simple root ${\alpha}_{v}\in B_{\TT_A}$ and $r\in W(R_{A})$.
Therefore we only need to show that $\widetilde{\alpha}_{v}+n\delta \in \Delta_{re}(\widetilde{R}_A)$ for all $v\in T_A$.
As in the following five cases, we can show that $\widetilde{\alpha}_{v}+n\delta \in \Delta_{re}(\widetilde{R}_A)$ for all $v\in T_A$;
\begin{enumerate}
\item If $v=(i,1)$, we have $\widetilde{\tau}^{n}_{1}(\widetilde{\alpha}_{v})=\widetilde{\alpha}_{v}+n\delta$.
\item If $v=(i,j)$ for $j\ge 2$, we have $\widetilde{\tau}^{n}_{(i,j-1)}(\widetilde{\alpha}_{v})=\widetilde{\alpha}_{v}+n\delta$.
\item If $v=1$ and $n=1$, we have $\widetilde{\alpha}_{1^{*}}=\widetilde{\alpha}_{1}+\delta$.
\item If $v=1$ and $n=2k+1$, we have $\widetilde{\tau}^{2k}_{1}(\widetilde{\alpha}_{1^{*}})=\widetilde{\alpha}_{1}+(2k+1)\delta$.
\item If $v=1$ and $n=2k$, we have $\widetilde{\tau}^{2k}_{1}(\widetilde{\alpha}_{1})=\widetilde{\alpha}_{1}+2k\delta$.
\end{enumerate}
Hence this lemma holds.
\qed
\end{pf}
Therefore we have finished the proof of this proposition. 
\qed 
\end{pf}

\section{Weyl groups as generalized Coxeter groups}
In this section, we express the Weyl group $W(\widetilde{R}_A)$ as a generalized Coxeter group. 
First we note the following fact.
\begin{prop}\label{prop:41}
Define a group $W({T}_A)$ by the following generators and the Coxeter relations attached to the diagram ${T}_A:$
\begin{description}
\item[{\bf Generators}] $\{{w}_{v}~|~v\in {T}_A\}$
\item[{\bf Relations}]
\begin{subequations}
\begin{equation}
{w}_{v}^2=1\quad \text{for all}\quad v\in {T}_A,
\end{equation}
\begin{equation}
w_{v}w_{v'}=w_{v'}w_{v}\quad \text{if}\quad I_{ {R}_A} ({\alpha}_{v}, {\alpha}_{v'})=0, 
\end{equation}
\begin{equation}
w_{v}w_{v'}w_{v}=w_{v'}w_{v}w_{v'}\quad \text{if}\quad I_{{R}_A}({\alpha}_{v}, {\alpha}_{v'})= -1.
\end{equation}
\end{subequations}
\end{description}
Then the correspondence ${w}_{v}\mapsto {r}_v$ for $v\in {T}_A$ induces an isomorphism of groups
\begin{equation}
W({T}_A)\cong {W}({R}_A).
\end{equation}
\end{prop}
\begin{pf}
It follows from the general result for Weyl groups associated to generalized Cartan matrices (cf. Proposition~3.13 in \cite{Kac1}).
\qed
\end{pf}
\begin{prop}\label{prop:42}
Define a group $W({T}_A)\ltimes K_0(R_A)$ by the following generators and the relations$:$
\begin{description}
\item[{\bf Generators}] $\{w_v, \tau_v~|~v\in T_A\}$
\item[{\bf Relations}] 
\begin{subequations}
\begin{equation}\label{4.3a}
w_{v}^2=1\quad \text{for all}\quad v\in T_A,
\end{equation}
\begin{equation}\label{4.3b}
w_{v}w_{v'}=w_{v'}w_{v}\quad \text{if}\quad I_{R_A}(\alpha_{v}, \alpha_{v'})= 0,
\end{equation}
\begin{equation}\label{4.3c}
w_{v}w_{v'}w_{v}=w_{v'}w_{v}w_{v'}\quad \text{if}\quad I_{R_A}(\alpha_{v}, \alpha_{v'})=-1,
\end{equation}
\begin{equation}\label{4.3d}
\tau_{v}\tau_{v'}=\tau_{v'}\tau_{v}\quad \text{for all}\quad v,v'\in T_A,
\end{equation}
\begin{equation}\label{4.3e}
w_{v}\tau_{v}w_{v}=\tau_{v}^{-1}\quad \text{for all}\quad v\in T_A,
\end{equation}
\begin{equation}\label{4.3f}
w_{v}\tau_{v'}=\tau_{v'}w_{v}\quad \text{if}\quad I_{R_A}(\alpha_{v},\alpha_{v'})= 0,
\end{equation}
\begin{equation}\label{4.3g}
w_{v}\tau_{v'}w_{v}=\tau_{v'}\tau_{v}\quad \text{if}\quad I_{R_A}(\alpha_{v},\alpha_{v'})= -1.
\end{equation}
\end{subequations}
\end{description}
Identify the subgroup generated by $\tau_v$, $v\in T_A$ with a free abelian group $K_0(R_A)$ expressed in multiplicative notation. 
\begin{enumerate}
\item
The correspondence ${w}_{v}\mapsto {r}_v$, $\tau_v\mapsto {\tau}_v$ for $v\in {T}_A$ induces an isomorphism 
of groups 
\begin{equation}
W({T}_A)\ltimes K_0(R_A)\cong W({R}_A)\ltimes K_0(R_A),
\end{equation}
where the semi-direct product in the right hand side is given by the natural inclusion $W(R_A)\hookrightarrow {\rm Aut}(K_0(R_A),I_{R_A})$.
\item
The correspondence ${w}_{v}\mapsto {r}_v$, $\tau_v\mapsto \widetilde{\tau}_v$ for $v\in {T}_A$ induces a surjective group homomorphism
\begin{equation}
W({T}_A)\ltimes K_0(R_A)\twoheadrightarrow W(\widetilde{R}_A),
\end{equation}
whose kernel is isomorphic to ${\rm rad}(I_{R_A})$.
\end{enumerate}
\end{prop}
\begin{pf}
The first statement is almost obvious since
we have ${r}_{v}({\alpha}_{v})=-{\alpha}_{v}$, ${r}_{v}({\alpha}_{v'})={\alpha}_{v'}$ if $I_{{R}_A}({\alpha}_{v}, {\alpha}_{v'}) = 0$ and 
${r}_{v}({\alpha}_{v'})={\alpha}_{v'}+{\alpha}_{v}$ if $I_{{R}_A}({\alpha}_{v}, {\alpha}_{v'}) = -1$.
The second statement follows from Proposition~\ref{prop:41} and Proposition~\ref{prop:34} that we have the surjective group homomorphism.
Theorem~\ref{relation W(R_A) and W(widetilde{R}_A)} and Proposition~\ref{translation} imply that 
the kernel is isomorphic to ${\rm rad}(I_{R_A})$.
\qed
\end{pf}
\begin{defn}\label{defn:43}
Define a group $W(\widetilde{T}_A)$ by the following generators and the generalized Coxeter relations attached to the diagram $\widetilde{T}_A:$
\begin{description}
\item[{\bf Generators}] $\{\widetilde{w}_{v}~|~v\in\widetilde{T}_A\}$
\item[{\bf Relations}]
\begin{subequations}
\begin{equation}\label{W0}\tag{{\bf W0}}
\widetilde{w}_{v}^2=1\quad \text{for all}\quad v\in\widetilde{T}_A,
\end{equation}
\begin{equation}\label{W1.0}\tag{{\bf W1.0}}
\widetilde{w}_{v}\widetilde{w}_{v'}=\widetilde{w}_{v'}\widetilde{w}_{v}\quad \text{if}\quad I_{\widetilde{R}_A} (\widetilde{\alpha}_{v},\widetilde{\alpha}_{v'})=0, 
\end{equation}
\begin{equation}\label{W1.1}\tag{{\bf W1.1}}
\widetilde{w}_{v}\widetilde{w}_{v'}\widetilde{w}_{v}=\widetilde{w}_{v'}\widetilde{w}_{v}\widetilde{w}_{v'}\quad \text{if}\quad I_{\widetilde{R}_A}(\widetilde{\alpha}_{v},\widetilde{\alpha}_{v'})= -1,
\end{equation}
\begin{equation}\label{W2}\tag{{\bf W2}}
\widetilde{w}_{(i,1)}\sigma_{1} \widetilde{w}_{(i,1)}\sigma_{1}=\sigma_{1}\widetilde{w}_{(i,1)}\sigma_{1}\widetilde{w}_{(i,1)},
\end{equation}
\begin{equation}\label{W3}\tag{{\bf W3}}
\begin{cases}
\widetilde{w}_{(i,1)}\sigma_{(j,1)}=\sigma_{(j,1)}\widetilde{w}_{(i,1)}\\
\widetilde{w}_{(j,1)}\sigma_{(i,1)}=\sigma_{(i,1)}\widetilde{w}_{(j,1)}\end{cases}
\quad \text{for all}\quad 1\le i<j\le r,
\end{equation}
\end{subequations}
where $\sigma_{1}:=\widetilde{w}_{1}\widetilde{w}_{1^*}$ and $\sigma_{(i,1)}:=\widetilde{w}_{(i,1)}\sigma_{1}\widetilde{w}_{(i,1)}\sigma_{1}^{-1}$ for all $i=1,\dots, r$.
\end{description}
\end{defn}
The conditions \eqref{W2} and \eqref{W3} are different from the definition in \cite{st:1}.
However we can deduce the original ones from \eqref{W2} and \eqref{W3}
under the conditon \eqref{W0}:
\begin{prop}[cf. Lemma~4.1 and Lemma~4.2 in \cite{yamada}]\label{lem 4.1 4.2 ver1}
Under the relation \ref{W0}, we have the following equivalences of relations$:$
\begin{subequations}
\begin{equation}
\widetilde{w}_{(i,1)}\sigma_{1} \widetilde{w}_{(i,1)}\sigma_{1}=\sigma_{1}\widetilde{w}_{(i,1)}\sigma_{1}\widetilde{w}_{(i,1)}
\quad \Longleftrightarrow\quad (\widetilde{w}_{1}\widetilde{w}_{(i,1)}\widetilde{w}_{1^{*}}\widetilde{w}_{(i,1)})^{3}=1,
\end{equation}
\begin{equation}
\widetilde{w}_{(i,1)}\sigma_{(j,1)}=\sigma_{(j,1)}\widetilde{w}_{(i,1)}
\quad \Longleftrightarrow\quad
(\widetilde{w}_{(i,1)}\widetilde{w}_{1}\widetilde{w}_{(i,1)}\widetilde{w}_{1^{*}}\widetilde{w}_{(j,1)}\widetilde{w}_{1^{*}})^{2}=1,
\end{equation}
\begin{equation}
\widetilde{w}_{(j,1)}\sigma_{(i,1)}=\sigma_{(i,1)}\widetilde{w}_{(j,1)}
\quad \Longleftrightarrow\quad 
(\widetilde{w}_{(i,1)}\widetilde{w}_{1^{*}}\widetilde{w}_{(i,1)}\widetilde{w}_{1}\widetilde{w}_{(j,1)}\widetilde{w}_{1})^{2}=1.
\end{equation}
\end{subequations}
\end{prop}
\begin{pf}
We have
\begin{eqnarray*}
& & (\widetilde{w}_{1}\widetilde{w}_{(i,1)}\widetilde{w}_{1^{*}}\widetilde{w}_{(i,1)})^{3}=1\\
&\Longleftrightarrow& \widetilde{w}_{1}\widetilde{w}_{(i,1)}\widetilde{w}_{1^{*}}\widetilde{w}_{(i,1)}\widetilde{w}_{1}\widetilde{w}_{(i,1)}\widetilde{w}_{1^{*}}\widetilde{w}_{(i,1)}\widetilde{w}_{1}\widetilde{w}_{(i,1)}\widetilde{w}_{1^{*}}\widetilde{w}_{(i,1)}=1\\
&\Longleftrightarrow& \widetilde{w}_{1}\widetilde{w}_{(i,1)}\widetilde{w}_{1^{*}}\widetilde{w}_{1}\widetilde{w}_{(i,1)}\widetilde{w}_{1}\widetilde{w}_{1^{*}}\widetilde{w}_{(i,1)}\widetilde{w}_{1}\widetilde{w}_{1^{*}}\widetilde{w}_{(i,1)}\widetilde{w}_{1^{*}}=1\quad (\text{by \eqref{W1.1}})\\
&\Longleftrightarrow& \widetilde{w}_{1}\widetilde{w}_{(i,1)}\sigma_{1}^{-1}\widetilde{w}_{(i,1)}\sigma_{1}\widetilde{w}_{(i,1)}\sigma_{1}\widetilde{w}_{(i,1)}\widetilde{w}_{1^{*}}=1\quad (\text{by the definition of }\sigma_{1})\\
&\Longleftrightarrow& \widetilde{w}_{(i,1)}\sigma_{1} \widetilde{w}_{(i,1)}\sigma_{1}=\sigma_{1}\widetilde{w}_{(i,1)}\widetilde{w}_{1}\widetilde{w}_{1^{*}}\widetilde{w}_{(i,1)}=\sigma_{1}\widetilde{w}_{(i,1)}\sigma_{1}\widetilde{w}_{(i,1)}.
\end{eqnarray*}
We have
\begin{eqnarray*}
& & \widetilde{w}_{(i,1)}\sigma_{(j,1)}=\sigma_{(j,1)}\widetilde{w}_{(i,1)}\\
&\Longleftrightarrow& \widetilde{w}_{(i,1)}\widetilde{w}_{(j,1)}\widetilde{w}_{1}\widetilde{w}_{1^{*}}\widetilde{w}_{(j,1)}\widetilde{w}_{1^{*}}\widetilde{w}_{1}=\widetilde{w}_{(j,1)}\widetilde{w}_{1}\widetilde{w}_{1^{*}}\widetilde{w}_{(j,1)}\widetilde{w}_{1^{*}}\widetilde{w}_{1}\widetilde{w}_{(i,1)}\quad (\text{by the definition of }\sigma_{(j,1)})\\
&\Longleftrightarrow& \widetilde{w}_{(j,1)}\widetilde{w}_{{(i,1)}}\widetilde{w}_{1}\widetilde{w}_{1^{*}}\widetilde{w}_{(j,1)}\widetilde{w}_{1^{*}}\widetilde{w}_{1}=\widetilde{w}_{(j,1)}\widetilde{w}_{1}\widetilde{w}_{1^{*}}\widetilde{w}_{(j,1)}\widetilde{w}_{1^{*}}\widetilde{w}_{1}\widetilde{w}_{(i,1)}\quad (\text{by \eqref{W1.1}})\\
&\Longleftrightarrow& \widetilde{w}_{{(i,1)}}\widetilde{w}_{1}\widetilde{w}_{1^{*}}\widetilde{w}_{(j,1)}\widetilde{w}_{1^{*}}\widetilde{w}_{1}=\widetilde{w}_{1}\widetilde{w}_{1^{*}}\widetilde{w}_{(j,1)}\widetilde{w}_{1^{*}}\widetilde{w}_{1}\widetilde{w}_{(i,1)}\\
&\Longleftrightarrow& \widetilde{w}_{1}\widetilde{w}_{{(i,1)}}\widetilde{w}_{1}\widetilde{w}_{1^{*}}\widetilde{w}_{(j,1)}\widetilde{w}_{1^{*}}=\widetilde{w}_{1^{*}}\widetilde{w}_{(j,1)}\widetilde{w}_{1^{*}}\widetilde{w}_{1}\widetilde{w}_{(i,1)}\widetilde{w}_{1}\\
&\Longleftrightarrow& \widetilde{w}_{(i,1)}\widetilde{w}_{1}\widetilde{w}_{(i,1)}\widetilde{w}_{1^{*}}\widetilde{w}_{(j,1)}\widetilde{w}_{1^{*}}=\widetilde{w}_{1^{*}}\widetilde{w}_{(j,1)}\widetilde{w}_{1^{*}}\widetilde{w}_{(i,1)}\widetilde{w}_{1}\widetilde{w}_{(i,1)}\quad (\text{by \eqref{W1.1}})\\
&\Longleftrightarrow& (\widetilde{w}_{(i,1)}\widetilde{w}_{1}\widetilde{w}_{(i,1)}\widetilde{w}_{1^{*}}\widetilde{w}_{(j,1)}\widetilde{w}_{1^{*}})^{2}=1.
\end{eqnarray*}
Similarly, we have
\begin{eqnarray*}
& & \widetilde{w}_{(j,1)}\sigma_{(i,1)}=\sigma_{(i,1)}\widetilde{w}_{(j,1)}\\
&\Longleftrightarrow& \widetilde{w}_{(j,1)}\widetilde{w}_{{(i,1)}}\widetilde{w}_{1}\widetilde{w}_{1^{*}}\widetilde{w}_{{(i,1)}}\widetilde{w}_{1^{*}}\widetilde{w}_{1}=\widetilde{w}_{{(i,1)}}\widetilde{w}_{1}\widetilde{w}_{1^{*}}\widetilde{w}_{{(i,1)}}\widetilde{w}_{1^{*}}\widetilde{w}_{1}\widetilde{w}_{(j,1)}\quad (\text{by the definition of }\sigma_{(i,1)})\\
&\Longleftrightarrow& \widetilde{w}_{(i,1)}\widetilde{w}_{(j,1)}\widetilde{w}_{1}\widetilde{w}_{1^{*}}\widetilde{w}_{{(i,1)}}\widetilde{w}_{1^{*}}\widetilde{w}_{1}=\widetilde{w}_{{(i,1)}}\widetilde{w}_{1}\widetilde{w}_{1^{*}}\widetilde{w}_{{(i,1)}}\widetilde{w}_{1^{*}}\widetilde{w}_{1}\widetilde{w}_{(j,1)}\quad (\text{by \eqref{W1.1}})\\
&\Longleftrightarrow& \widetilde{w}_{(j,1)}\widetilde{w}_{1}\widetilde{w}_{1^{*}}\widetilde{w}_{{(i,1)}}\widetilde{w}_{1^{*}}\widetilde{w}_{1}=\widetilde{w}_{1}\widetilde{w}_{1^{*}}\widetilde{w}_{{(i,1)}}\widetilde{w}_{1^{*}}\widetilde{w}_{1}\widetilde{w}_{(j,1)}\\
&\Longleftrightarrow& \widetilde{w}_{1}\widetilde{w}_{(j,1)}\widetilde{w}_{1}\widetilde{w}_{1^{*}}\widetilde{w}_{{(i,1)}}\widetilde{w}_{1^{*}}=\widetilde{w}_{1^{*}}\widetilde{w}_{{(i,1)}}\widetilde{w}_{1^{*}}\widetilde{w}_{1}\widetilde{w}_{(j,1)}\widetilde{w}_{1}\\
&\Longleftrightarrow& \widetilde{w}_{1}\widetilde{w}_{(j,1)}\widetilde{w}_{1}\widetilde{w}_{(i,1)}\widetilde{w}_{1^{*}}\widetilde{w}_{(i,1)}=\widetilde{w}_{{(i,1)}}\widetilde{w}_{1^{*}}\widetilde{w}_{(i,1)}\widetilde{w}_{1}\widetilde{w}_{(j,1)}\widetilde{w}_{1}\quad (\text{by \eqref{W1.1}})\\
&\Longleftrightarrow& (\widetilde{w}_{(i,1)}\widetilde{w}_{1^{*}}\widetilde{w}_{(i,1)}\widetilde{w}_{1}\widetilde{w}_{(j,1)}\widetilde{w}_{1})^{2}=1.
\end{eqnarray*}
We have finished the proof of the proposition.
\qed
\end{pf}
Note that the Coxeter--Dynkin diagram $\widetilde{T}_A$ is symmetric under the permutation 
\begin{equation}
1^{*}\mapsto 1, \ 1\mapsto 1^{*}, \ v\mapsto v \quad \text{if} \quad v\neq 1, 1^{*}.
\end{equation}
This symmetry of $\widetilde{T}_A$ induces the automorphism on $W(\widetilde{R}_A)$ which sends $\sigma_{1}$ to $\sigma^{-1}_{1}$ 
and hence $W(\widetilde{T}_A)$ depends only on the Coxeter--Dynkin diagram $\widetilde{T}_A$.
\begin{thm}\label{main}
We have an isomorphism of groups
\begin{equation}
W(\widetilde{T}_A)\cong W(R_A)\ltimes K_0(R_A).
\end{equation}
In particular, $W(\widetilde{T}_A)\cong W(\widetilde{R}_A)$ if $\chi_A\ne 0$ and 
$W(\widetilde{T}_A)$ is isomorphic to the hyperbolic extension of the elliptic Weyl group $W(\widetilde{R}_A)$ if $\chi_A=0$.
\end{thm}
\begin{pf}
We shall show that $W(\widetilde{T}_A)\cong W(T_A)\ltimes K_0(R_A)$, which yields the statement by Proposition~\ref{prop:42}.
Define an element $\sigma_{v}\in W(\widetilde{T}_A)$ for each vertex $v\in T_A$ by induction as follows$:$
\begin{itemize}
\item
For verteces $1$ and $(i,1)$ for $i=1,\dots, r$, in Definition~\ref{defn:43} we have already set  
\begin{subequations}
\begin{equation}
\sigma_{1}:=\widetilde{w}_{1}\widetilde{w}_{1^{*}}.
\end{equation}
\begin{equation}
\sigma_{(i,1)}:=\widetilde{w}_{{(i,1)}}\sigma_{1}\widetilde{w}_{{(i,1)}}\sigma_{1}^{-1},
\quad i=1,\dots, r.
\end{equation}
\item
Set
\begin{equation}
\sigma_{(i,j)}:=\widetilde{w}_{{(i,j)}}\sigma_{(i,j-1)}\widetilde{w}_{{(i,j)}}\sigma_{(i,j-1)}^{-1},
\quad i=1,\dots, r,\ j=2,\dots, a_i-1.
\end{equation}
\end{subequations}
\end{itemize}
\begin{lem}\label{surjective:4.6}
The correspondences $w_{v}\mapsto \widetilde{w}_{v}, \tau_{v}\mapsto {\sigma}_{v}$ induce the surjective group homomorphism$:$
\begin{equation}
W(T_A)\ltimes K_0(R_A)\twoheadrightarrow W(\widetilde{T}_A).
\end{equation}
\end{lem}
\begin{pf}
We shall check that the elements in $W(\widetilde{T}_A)$ satisfy the relations of $W(T_A)\ltimes K_0(R_A)$.
First, we check the relation \eqref{4.3d}:
\begin{sublem}\label{S commute}
The elements $\sigma_{v}$ for $v\in T_A$ commute with each other.
\end{sublem}
\begin{pf}

First, we prove $\sigma_{1}$ and $\sigma_{(i,1)}$ are commuting.
\begin{eqnarray*}
\sigma_{(i,1)}\sigma_{1}
&=& \widetilde{w}_{(i,1)}\sigma_{1} \widetilde{w}_{(i,1)} \sigma_{1}^{-1} \sigma_{1} \quad (\text{by the definition of } \sigma_{(i,1)})\\
&=& \widetilde{w}_{(i,1)}\sigma_{1} \widetilde{w}_{(i,1)}\\
&=& \sigma_{1} \widetilde{w}_{(i,1)} \sigma_{1} \widetilde{w}_{(i,1)} \sigma_{1}^{-1} \\
&=& \sigma_{1}\sigma_{(i,1)} \quad (\text{by the definition of } \sigma_{(i,1)}).
\end{eqnarray*}
Next, we prove that $\widetilde{w}_{(i,2)}$ and $\sigma_{(i,1)}$ satisfy the following relation:
\begin{equation*}
\widetilde{w}_{(i,2)}\sigma_{(i,1)}\widetilde{w}_{(i,2)}\sigma_{(i,1)} = \sigma_{(i,1)}\widetilde{w}_{(i,2)}\sigma_{(i,1)}\widetilde{w}_{(i,2)}.
\end{equation*}
By \eqref{W1.0} and commutativity of $\sigma_{1}$ and $\sigma_{(i,1)}$, one obtains the following:
\begin{eqnarray*}
& & \widetilde{w}_{(i,2)}\sigma_{(i,1)} \widetilde{w}_{(i,2)} \sigma_{(i,1)}\\
&=&\widetilde{w}_{(i,2)}\widetilde{w}_{(i,1)} \sigma_{1} \widetilde{w}_{(i,1)} \sigma_{1}^{-1} \widetilde{w}_{(i,2)} \sigma_{(i,1)} \quad 
\text{(Definition of $\sigma_{(i,1)}$)}\\
&=&\widetilde{w}_{(i,2)}\widetilde{w}_{(i,1)} \sigma_{1} \widetilde{w}_{(i,1)} \widetilde{w}_{(i,2)} \sigma_{(i,1)}\sigma_{1}^{-1} \quad 
\text{(Commutativity of $\sigma_{1}$ and $\sigma_{(i,1)}$)}\\
&=&(\widetilde{w}_{(i,2)}\widetilde{w}_{(i,1)} \sigma_{1} \widetilde{w}_{(i,1)} \widetilde{w}_{(i,2)} \widetilde{w}_{(i,1)} \sigma_{1}  \widetilde{w}_{(i,1)} \sigma_{1}^{-1} )\sigma_{1}^{-1} \quad \text{(Definition of $\sigma_{(i,1)}$)}
\\ 
&=&(\widetilde{w}_{(i,2)}\widetilde{w}_{(i,1)} \sigma_{1} \widetilde{w}_{(i,2)} \widetilde{w}_{(i,1)} \widetilde{w}_{(i,2)} \sigma_{1}  \widetilde{w}_{(i,1)})\sigma_{1}^{-2} \quad \text{(By the relation \eqref{W1.1})}
\\
&=&(\widetilde{w}_{(i,2)}\widetilde{w}_{(i,1)}\widetilde{w}_{(i,2)} \sigma_{1} \widetilde{w}_{(i,1)} \sigma_{1} \widetilde{w}_{(i,2)} \widetilde{w}_{(i,1)})\sigma_{1}^{-2} \\
&=&(\widetilde{w}_{(i,1)}\widetilde{w}_{(i,2)}\widetilde{w}_{(i,1)} \sigma_{1} \widetilde{w}_{(i,1)} \sigma_{1} \widetilde{w}_{(i,2)} \widetilde{w}_{(i,1)})\sigma_{1}^{-2} \quad \text{(By the relation \eqref{W1.1})}
\\
&=& (\widetilde{w}_{(i,1)} \widetilde{w}_{(i,2)} \sigma_{1}\widetilde{w}_{(i,1)}\sigma_{1} \widetilde{w}_{(i,1)} \widetilde{w}_{(i,2)} \widetilde{w}_{(i,1)}) \sigma_{1}^{-2} \quad \text{(Commutativity of $\sigma_{1}$ and $\sigma_{(i,1)}$)}
\\
&=& (\widetilde{w}_{(i,1)} \widetilde{w}_{(i,2)} \sigma_{1}\widetilde{w}_{(i,1)}\sigma_{1} \widetilde{w}_{(i,2)} \widetilde{w}_{(i,1)} \widetilde{w}_{(i,2)}) \sigma_{1}^{-2} \quad \text{(By the relation \eqref{W1.1})}
\\
&=& (\widetilde{w}_{(i,1)} \sigma_{1} \widetilde{w}_{(i,2)}  \widetilde{w}_{(i,1)}\widetilde{w}_{(i,2)} \sigma_{1} \widetilde{w}_{(i,1)} \widetilde{w}_{(i,2)}) \sigma_{1}^{-2}\\
&=& (\widetilde{w}_{(i,1)} \sigma_{1} \widetilde{w}_{(i,1)}  \widetilde{w}_{(i,2)}\widetilde{w}_{(i,1)} \sigma_{1} \widetilde{w}_{(i,1)} \widetilde{w}_{(i,2)}) \sigma_{1}^{-2} \quad \text{(By the relation \eqref{W1.1})}
\\
&=& (\widetilde{w}_{(i,1)} \sigma_{1} \widetilde{w}_{(i,1)}  \widetilde{w}_{(i,2)}\widetilde{w}_{(i,1)} \sigma_{1} \widetilde{w}_{(i,1)} \sigma_{1}^{-1} \widetilde{w}_{(i,2)}) \sigma_{1}^{-1}\\
&=& (\widetilde{w}_{(i,1)} \sigma_{1} \widetilde{w}_{(i,1)}  \widetilde{w}_{(i,2)} \sigma_{(i,1)} \widetilde{w}_{(i,2)}) \sigma_{1}^{-1}
\quad \text{(Definition of $\sigma_{(i,1)}$)}
\\
&=& \widetilde{w}_{(i,1)} \sigma_{1} \widetilde{w}_{(i,1)} \sigma_{1}^{-1} \widetilde{w}_{(i,2)} \sigma_{(i,1)} \widetilde{w}_{(i,2)}
\quad \text{(Commutativity of $\sigma_{1}$ and $\sigma_{(i,1)}$)}\\
&=& \sigma_{(i,1)} \widetilde{w}_{(i,2)} \sigma_{(i,1)} \widetilde{w}_{(i,2)} \quad \text{(Definition of $\sigma_{(i,1)}$)}.
\end{eqnarray*}
From this relation and the definition of $\sigma_{(i,2)}$, one has
\begin{equation*}
\sigma_{(i,2)}\sigma_{(i,1)} = \widetilde{w}_{(i,2)} \sigma_{(i,1)} \widetilde{w}_{(i,2)} = \sigma_{(i,1)} \widetilde{w}_{(i,2)} \sigma_{(i,1)} \widetilde{w}_{(i,2)} \sigma_{(i,1)}^{-1} = \sigma_{(i,1)} \sigma_{(i,2)}.
\end{equation*}
Hence we have
\begin{equation*}
\sigma_{(i,1)}\sigma_{(i,2)}=\sigma_{(i,2)}\sigma_{(i,1)}.
\end{equation*}
By similar calculation, we obtain
\begin{equation*}
\widetilde{w}_{(i,j+1)} \sigma_{(i,j)} \widetilde{w}_{(i,j+1)} \sigma_{(i,j)} = \sigma_{(i,j)} \widetilde{w}_{(i,j+1)} \sigma_{(i,j)} \widetilde{w}_{(i,j+1)}\quad \text{for}\quad j=2,\cdots,a_{i}-1,
\end{equation*}
and $\sigma_{(i,j)}\sigma_{(i,j+1)}=\sigma_{(i,j+1)}\sigma_{(i,j)}$.

By the definition of $\sigma_{(i,j)}$, it is written by products of 
$\widetilde{w}_{1}, \widetilde{w}_{1^{*}}, \widetilde{w}_{(i,k)}$ for $k=1,\dots j$. 
Then $\sigma_{(i,j)}$ commutes with $\widetilde{w}_{(i,j+2)}$ by \eqref{W1.0}. 
Then we have $\sigma_{(i,j+2)}\sigma_{(i,j)}=\sigma_{(i,j)}\sigma_{(i,j+2)}$. 
Therefore, by same calculations,
we obtain the commutativity of $\sigma_{(i,j)}$ and $\sigma_{(i,k)}$ for 
$j,k=1,2,\cdots, a_{i}-1$.
\qed
\end{pf}

Next we ckeck that the relations \eqref{4.3e}, \eqref{4.3f} and \eqref{4.3g}:
\begin{sublem}
For $v,v'\in T_{A}$, we have
\begin{subequations}
\begin{equation}\label{adjunction v'=v}
\widetilde{w}_{v}\sigma_{v}\widetilde{w}_{v}=\sigma_{v}^{-1},
\end{equation}
\begin{equation}\label{adjunction 0}
\widetilde{w}_{v}\sigma_{v'}\widetilde{w}_{v}=\sigma_{v'}\quad \text{if}\quad I_{\widetilde{R}_A}(\widetilde{\alpha}_{v},\widetilde{\alpha}_{v'})=I_{R_A}({\alpha}_{v},{\alpha}_{v'})=0,
\end{equation}
\begin{equation}\label{adjunction -1}
\widetilde{w}_{v}\sigma_{v'}\widetilde{w}_{v}=\sigma_v \sigma_{v'}\quad \text{if}\quad I_{\widetilde{R}_A}(\widetilde{\alpha}_{v},\widetilde{\alpha}_{v'})=I_{R_A}({\alpha}_{v},{\alpha}_{v'})=-1.
\end{equation}
\end{subequations}
\end{sublem}
\begin{pf}
First, we show the relation \eqref{adjunction v'=v}. 
For the case that $v=1$, we have
\begin{equation*}
\widetilde{w}_{1}\sigma_{1}\widetilde{w}_{1}=\widetilde{w}_{1}\widetilde{w}_{1}\widetilde{w}_{1^{*}}\widetilde{w}_{1}=(\widetilde{w}_{1}\widetilde{w}_{1^{*}})^{-1}=\sigma^{-1}_{1}.
\end{equation*}
For the case that $v=(i,j)$, we have
\begin{eqnarray*}
&&\widetilde{w}_{(i,j)}\sigma_{(i,j)}\widetilde{w}_{(i,j)}\\
&=&\widetilde{w}_{(i,j)} (\widetilde{w}_{(i,j)}\sigma_{(i,j-1)} \widetilde{w}_{(i,j)} \sigma^{-1}_{(i,j-1)})\widetilde{w}_{(i,j)}\\
&=&(\widetilde{w}_{(i,j)}\sigma_{(i,j-1)} \widetilde{w}_{(i,j)} \sigma^{-1}_{(i,j-1)})^{-1}=\sigma_{(i,j)}^{-1}.
\end{eqnarray*}
Here we write $\sigma_{(i,0)}:=\sigma_{1}$.

Second, we shall show the relation \eqref{adjunction 0}.
By the relation \eqref{W3}, we have 
\begin{equation*}
\widetilde{w}_{(i,1)} \sigma_{(i',1)} =\sigma_{(i',1)} \widetilde{w}_{(i,1)}.\quad 
\end{equation*}
Assume the following equation:
\begin{equation*}
\widetilde{w}_{(i,j)} \sigma_{(i',k)} = \sigma_{(i',k)} \widetilde{w}_{(i,j)}.
\end{equation*}
Under this assumption, we have
\begin{eqnarray*}
\widetilde{w}_{(i,j)}\sigma_{(i',k+1)}\widetilde{w}_{(i,j)}
&=& \widetilde{w}_{(i,j)} ( \widetilde{w}_{(i',k+1)}\sigma_{(i',k)} \widetilde{w}_{(i',k+1)} \sigma^{-1}_{(i',k)})\widetilde{w}_{(i,j)}\\
&=& ( \widetilde{w}_{(i',k+1)}\sigma_{(i',k)} \widetilde{w}_{(i',k+1)} \sigma^{-1}_{(i',k)}) \widetilde{w}_{(i,j)}\widetilde{w}_{(i,j)}\\
&=& \sigma_{(i',k+1)}.
\end{eqnarray*}
Finally, we shall show the relation \eqref{adjunction -1}.
We have
\begin{eqnarray*}
& & \widetilde{w}_{(i,j+1)}\sigma_{(i,j)}\widetilde{w}_{(i,j+1)}\\
&=&(\widetilde{w}_{(i,j+1)}\sigma_{(i,j)}\widetilde{w}_{(i,j+1)} \sigma_{(i,j)})\sigma^{-1}_{(i,j)} \\
&=&(\sigma_{(i,j)}\widetilde{w}_{(i,j+1)} \sigma_{(i,j)}  \widetilde{w}_{(i,j+1)})\sigma^{-1}_{(i,j)}\quad \text{(By Lemma \ref{S commute})}\\
&=&\sigma_{(i,j)}(\widetilde{w}_{(i,j+1)} \sigma_{(i,j)}  \widetilde{w}_{(i,j+1)}\sigma^{-1}_{(i,j)})\quad \text{(Definition of $\sigma_{(i,j+1)}$)}\\
&=&\sigma_{(i,j)}\sigma_{(i,j+1)}=\sigma_{(i,j+1)}\sigma_{(i,j)}.  
\end{eqnarray*}
We have finished the proof of the sublemma.
\qed
\end{pf}
We have finished the proof of Lemma \ref{surjective:4.6}.
\qed
\end{pf}

To construct the inverse homomorphism, let us define the
element $w_{1^{*}}$ of the group $W(T_A)\ltimes K_{0}(R_A)$ by
\begin{equation}
w_{1^{*}}:=w_{1}^{-1}\tau_{1}=w_{1}\tau_{1}.
\end{equation}

\begin{lem}\label{lem:inverse}
The elements $w_{v}$, $v\in T_A$ and $w_{1^{*}}$ of $W({T}_A)\ltimes K_0(R_A)$
satisfy the relations in Definition \ref{defn:43}. 
\end{lem}

\begin{pf}
First, we show the relation \eqref{W0}.
We only have to check that $w_{1^{*}}^{2}=1$. By the relation \eqref{4.3e}, we have 
\begin{equation*}
w_{1^{*}}^{2}=w_{1}\tau_{1}w_{1}\tau_{1}=1.
\end{equation*}

Second we show the relation \eqref{W1.0}.
We only have to check that $w_{1^{*}}w_{(i,j)}=w_{(i,j)}w_{1^{*}}$ for $j=2,\dots, a_{i}-1$.
By the relation \eqref{4.3f}, we have 
\begin{equation*}
w_{1^{*}}w_{(i,j)}=w_{1}\tau_{1}w_{(i,j)}=w_{1}w_{(i,j)}\tau_{1}=w_{(i,j)}w_{1}\tau_{1}
=w_{(i,j)}w_{1^{*}}.
\end{equation*}

Third we show the relation \eqref{W1.1}.
We only have to show the relation
\begin{equation*}
w_{(i,1)} w_{1^{*}} w_{(i,1)} = w_{1^{*}} w_{(i,1)} w_{1^{*}},
\end{equation*}
which follows from the calculation
\begin{eqnarray*}
w_{1^{*}} w_{(i,1)} w_{1^{*}}
&=& w_{1}w_{(i,1)}(w_{(i,1)} w_{1} w_{1^{*}} w_{(i,1)})w_{1^{*}}\\
&=& w_{1}w_{(i,1)}(w_{(i,1)} \tau_{1} w_{(i,1)} \tau_{1}^{-1} \tau_{1})w_{1^{*}}\\
&=& w_{1}w_{(i,1)}(\tau_{(i,1)} \tau_{1} )w_{1^{*}}\\
&=&w_{1}w_{(i,1)}(w_{1} \tau_{(i,1)} w_{1})w_{1^{*}}\\
&=& w_{1}w_{(i,1)}(w_{1} w_{(i,1)} \tau_{1} w_{(i,1)} \tau_{1}^{-1} w_{1})w_{1^{*}} \\
&=& w_{1}w_{(i,1)}(w_{1} w_{(i,1)} w_{1} w_{1^{*}} w_{(i,1)} (w_{1^{*}})^{-1} w_{1}^{-1} w_{1})w_{1^{*}} \\
&=& w_{1}w_{(i,1)}(w_{(i,1)} w_{1} w_{(i,1)} w_{1^{*}} w_{(i,1)} (w_{1^{*}})^{-1})w_{1^{*}}\\
&=& w_{(i,1)} w_{1^{*}}w_{(i,1)}.
\end{eqnarray*}

Finally, we show the relation \eqref{W2}.
Since the elements $w_{(i,1)}$ and $\tau_{1}$ satisfy the relation
\begin{equation*}
w_{(i,1)} \tau_{1} w_{(i,1)} = \tau_{(i,1)} \tau_{1},
\end{equation*}
by the commutativity of $\tau_{(i,1)}$ and $\tau_{1}$, we have
\begin{equation*}
w_{(i,1)}\tau_{1}w_{(i,1)}\tau_{1}=\tau_{(i,1)} \tau_{1} \tau_{1} =\tau_{1} \tau_{(i,1)} \tau_{1} = \tau_{1} w_{(i,1)} \tau_{1} w_{(i,1)}.
\end{equation*}
The relation \eqref{W3} is clearly satisfied since it is a part of the relation~\eqref{4.3f}.
We have finished the proof of the lemma.
\qed
\end{pf}
Thus we have completed the proof of Theorem~\ref{main}.
\qed
\end{pf}
%

\section{Cuspidal Artin Groups}

In this section, we obtain a relation between the generalized Coxeter group $W(\widetilde{T}_A)$ and the fundamental group
of regular orbit space for $W(R_A)\ltimes K_{0}(R_A)$.
\begin{defn}\label{def G'}
Define a group $G(\widetilde{T}_A)$ by the following generators and the generalized Coxeter relations attached to the diagram $\widetilde{T}_A:$
\begin{description}
\item[{\bf Generators}] $\{\widetilde{g}_v~|~v\in\widetilde{T}_A\}$
\item[{\bf Relations}] 
\begin{subequations}
\begin{equation}\label{A1.0}\tag{{\bf A1.0}}
\widetilde{g}_{v}\widetilde{g}_{v'}=\widetilde{g}_{v'}\widetilde{g}_{v}\quad \text{if}\quad I_{\widetilde{R}_A}(\widetilde{\alpha}_{v},\widetilde{\alpha}_{v'})= 0,
\end{equation}
\begin{equation}\label{A1.1}\tag{{\bf A1.1}}
\widetilde{g}_{v}\widetilde{g}_{v'}\widetilde{g}_{v}=\widetilde{g}_{v'}\widetilde{g}_{v}\widetilde{g}_{v'}\quad \text{if}\quad I_{\widetilde{R}_A}(\widetilde{\alpha}_{v},\widetilde{\alpha}_{v'})=-1,
\end{equation}
\begin{equation}\label{A2}\tag{{\bf A2}}
\widetilde{g}_{(i,1)}\widetilde{\rho}_{1} \widetilde{g}_{(i,1)}\widetilde{\rho}_{1}=\widetilde{\rho}_{1}\widetilde{g}_{(i,1)}\widetilde{\rho}_{1}\widetilde{g}_{(i,1)}\quad \text{for all } i=1,\dots, r,
\end{equation}
\begin{equation}\label{A3}\tag{{\bf A3}}
\begin{cases}
\widetilde{g}_{(i,1)}\widetilde{\rho}_{(j,1)}=\widetilde{\rho}_{(j,1)}\widetilde{g}_{(i,1)}\\
\widetilde{g}_{(j,1)}\widetilde{\rho}_{(i,1)}=\widetilde{\rho}_{(i,1)}\widetilde{g}_{(j,1)}
\end{cases}
\quad \text{for all}\quad 1\le i<j\le r.
\end{equation}
\end{subequations}
where $\widetilde{\rho}_{1}:=\widetilde{g}_{1}\widetilde{g}_{1^*}$ and $\widetilde{\rho}_{(i,1)}:=\widetilde{g}_{(i,1)}\widetilde{\rho}_{1}\widetilde{g}_{(i,1)}\widetilde{\rho}_{1}^{-1}$ for all $i=1,\dots, r$.
\end{description}
\end{defn}
\begin{defn}\label{defn:cuspidal Artin}
Let the notations be as above.
\begin{enumerate}
\item
If $\chi_A=0$, then the group $G(\widetilde{T}_A)$ is called the {\it elliptic Artin group} of type $A$.
\item
If $\chi_A<0$, then the group $G(\widetilde{T}_A)$ is called the {\it cuspidal Artin group} of type $A$.
\end{enumerate}
\end{defn}
\begin{rem}
It turns out later that the group $G(\widetilde{T}_A)$ is an affine Artin group by Theorem~\ref{G isom G'} if $\chi_A>0$.
\end{rem}
In addition to $\widetilde{\rho}_{1}$, $\widetilde{\rho}_{(i,1)}$, we also define the element $\widetilde{\rho}_{(i,j+1)}$ inductively as follows:
\begin{equation}
\widetilde{\rho}_{(i,j+1)}:=\widetilde{g}_{(i,j+1)}\widetilde{\rho}_{(i,j)}\widetilde{g}_{(i,j+1)}\widetilde{\rho}_{(i,j)}^{-1}, \quad 
i=1,\dots, r,\ 
j=1,\dots, a_{i}-2.
\end{equation}
\begin{prop}\label{W' isom G' bar}
The correspondence $\widetilde{g}_v\mapsto \widetilde{w}_v$ for $v\in\widetilde{T}_A$ induces a surjective group homomorphism 
\begin{equation}
G(\widetilde{T}_A)\twoheadrightarrow W(\widetilde{T}_A), 
\end{equation}
which yields an isomorphism 
\begin{equation}
G(\widetilde{T}_A)\left/\langle \widetilde{g}_v^2~\vert~v\in\widetilde{T}_A\rangle\right.\cong  W(\widetilde{T}_A).
\end{equation}
\end{prop}
\begin{pf}
It is obvious from Definition~\ref{defn:43}.
\qed
\end{pf}
\begin{defn}\label{defn of Tits cone}
Define a complex manifold $\E(R_A)$ by
\begin{equation}
\E(R_A):=\{h\in K_0(R_A)_\CC^*~\vert~ {\rm Im}(h)\in C(R_A)\},
\end{equation}
where $C(R_A)$ is the topological interior of the {\it Tits cone} $\overline{C}(R_A)$ of $R_A:$
\begin{equation}
\overline{C}(R_A):=\bigcup_{w\in W(R_{A})}w \left(\{ h\in K_0(R_A)_\CC^*~\vert~ h(\alpha_v)\ge 0, \  \text{for all} \ v \in T_{A} \}\right).
\end{equation}
Set
\begin{equation}
\E(R_A)^{reg}:=\E(R_A)\setminus\bigcup_{\alpha\in\Delta_{re}(R_A),n\in\ZZ} H_{\alpha, n}.
\end{equation}
where we denote by $H_{\alpha, n}$ the reflection hyperplane associated to $\widetilde{T}_A$, i.e.,
\begin{equation}
H_{\alpha, n}:=\{h\in K_0(R_A)_\CC^*~|h(\alpha)=n\}.
\end{equation}
\end{defn}
The group $W(R_A)\ltimes K_0(R_A)$ naturally acts on $\E(R_A)$ 
in a proper discontinuous way. 
It is known that the action is free on $\E(R_A)^{reg}$.
\begin{defn}
Define a group $G(\widetilde{R}_A)$ as the fundamental group of the regular orbit space$:$ 
\begin{equation}
G(\widetilde{R}_A):= \pi_{1}(\E(R_A)^{reg} /( W(R_A)\ltimes K_0(R_A)),*).
\end{equation}
\end{defn}
\begin{rem}
Since the complex manifold $\E(R_A)^{reg}$ is connected,
the group $G(\widetilde{R}_A)$ does not depend on the base point $*$.
\end{rem}
By definition of fundamental groups, we have the following commutative diagram of groups$:$
{\small
\[
\xymatrix{
 & & \{1\}\ar[d] &  \{1\}\ar[d]&\\
\{1\}\ar[r] & \pi_{1}\left(\E(R_A)^{reg},*\right)\ar[r]\ar@{=}[d]& \pi_{1}\left(\E(R_A)^{reg}/K_0(R_A),*\right)\ar[r]\ar[d] & K_0(R_A)\ar[r]\ar[d] &\{1\}\\
\{1\}\ar[r] & \pi_{1}\left(\E(R_A)^{reg},*\right)\ar[r]& G(\widetilde{R}_A)\ar[r]\ar[d] &W(R_A)\ltimes K_0(R_A)\ar[r]\ar[d] &\{1\}\\
& & W(R_A)\ar@{=}[r]\ar[d]& W(R_A)\ar[d]&\\
& & \{1\}& \{1\}&
}
\]
}
Generalizing the result for $\chi_A=0$ by Yamada \cite{yamada}, we obtain the following$:$ 
\begin{thm}\label{G isom G'}
There exists an isomorphism of groups
\begin{equation}
G(\widetilde{T}_A)\cong G(\widetilde{R}_A).
\end{equation}
\end{thm}
\begin{pf}
We can prove exactly in the same way as Yamada \cite{yamada}. 
The key is the following description of $G(\widetilde{R}_A)$ by Van der Lek \cite{van}.
\begin{prop}
The group $G(\widetilde{R}_A)$ is described by the following generators and relations$:$
\begin{description}
\item[{\bf Generators}] $\{g_v, \rho_v~|~v\in T_A\}$
\item[{\bf Relations}] 
\begin{subequations}
\begin{equation}\label{E 1}
g_{v}g_{v'}=g_{v'}g_{v}\quad \text{if}\quad I_{R_A}(\alpha_{v}, \alpha_{v'})= 0,
\end{equation}
\begin{equation}\label{E 1-2}
g_{v}g_{v'}g_{v}=g_{v'}g_{v}g_{v'}\quad \text{if}\quad I_{R_A}(\alpha_{v}, \alpha_{v'})=-1,
\end{equation}
\begin{equation}\label{E c}
\rho_{v}\rho_{v'}=\rho_{v'}\rho_{v}\quad \text{for all}\quad v,v'\in T_A,
\end{equation}
\begin{equation}\label{E 3}
g_{v}\rho_{v'}=\rho_{v'}g_{v}\quad \text{if}\quad I_{R_A}(\alpha_{v},\alpha_{v'})= 0,
\end{equation}
\begin{equation}\label{E a}
g_{v}\rho_{v'}g_{v}=\rho_{v'}\rho_{v}\quad \text{if}\quad I_{R_A}(\alpha_{v},\alpha_{v'})= -1.
\end{equation}
\end{subequations}
\end{description}
\end{prop}
\begin{pf}
See Theorem in \cite{van}.
\qed
\end{pf}
\begin{cor}\label{surj:5.8}
The correspondences ${g}_v\mapsto {w}_v$, $\rho_v\mapsto \tau_v$ for $v\in {T}_A$ induces a surjective group homomorphism 
\begin{equation}
G(\widetilde{R}_A)\twoheadrightarrow W(T_A)\ltimes K_0(T_A), 
\end{equation}
which yields an isomorphism 
\begin{equation}
G(\widetilde{R}_A)\left/\langle {g}_v^2,\ g_v\rho_v g_v\rho_v~\vert~v\in\widetilde{T}_A\rangle\right.\cong W(R_A)\ltimes K_0(R_A).
\end{equation}
\end{cor}
\begin{pf}
It is obvious from Proposition~\ref{prop:42}.
\qed
\end{pf}
\begin{lem}
The correspondences $g_{v}\mapsto \widetilde{g}_{v}, \rho_{v}\mapsto \widetilde{\rho}_{v}$ induce the surjective group homomorphism$:$
\begin{equation}
G(\widetilde{R}_A)\twoheadrightarrow G(\widetilde{T}_A).
\end{equation}
\end{lem}

\begin{pf}
We can easily check that the elements $\widetilde{g}_{v}$ satisfy
the relations \eqref{E 1} and \eqref{E 1-2} by using the relations \eqref{A1.0} and \eqref{A1.1}.

First we show the relation \eqref{E a}. We have
\begin{eqnarray*}
\widetilde{g}_{1}\widetilde{\rho}_{(i,1)}\widetilde{g}_1
&=& \widetilde{g}_{1}\widetilde{g}_{(i,1)}\widetilde{\rho}_{1}\widetilde{g}_{(i,1)}\widetilde{\rho}_{1}^{-1}\widetilde{g}_1\quad (\text{by the definition of }\widetilde{\rho}_{(i,1)})\\
&=& \widetilde{g}_{1}\widetilde{g}_{(i,1)}\widetilde{g}_{1}\widetilde{g}_{1^*}\widetilde{g}_{(i,1)}\widetilde{g}_{1^*}^{-1}\quad (\text{by the definition of }\widetilde{\rho}_{1})\\
&=& \widetilde{g}_{(i,1)}\widetilde{g}_{1}\widetilde{g}_{(i,1)}\widetilde{g}_{1^*}\widetilde{g}_{(i,1)}\widetilde{g}_{1^*}^{-1}\quad (\text{by \eqref{A1.1}})\\
&=& \widetilde{g}_{(i,1)}\widetilde{g}_{1}\widetilde{g}_{1^*}\widetilde{g}_{(i,1)}\quad (\text{by \eqref{A1.1}})\\
&=& \widetilde{g}_{(i,1)}\widetilde{\rho}_{1}\widetilde{g}_{(i,1)}\quad (\text{by the definition of }\widetilde{\rho}_{1})\\
&=& \widetilde{\rho}_{(i,1)}\widetilde{\rho}_{1}\quad (\text{by the definition of }\widetilde{\rho}_{(i,1)}).
\end{eqnarray*}

We can easily check that $\widetilde{g}_{v}$ and $\widetilde{\rho}_{v}$ satisfy
the relation \eqref{E 3} by the relation \eqref{A3}.

Second we show the relation \eqref{E c}.
By the relation \eqref{E a} and the definition of $\widetilde{\rho}_{v}$ for $v\in T_{A}$, we have
\begin{eqnarray*}
& & \widetilde{g}_{(i,2)}\widetilde{\rho}_{(i,1)} \widetilde{g}_{(i,2)} \widetilde{\rho}_{(i,1)}\\
&=&\widetilde{g}_{(i,2)}\widetilde{g}_{(i,1)} \widetilde{\rho}_{1} \widetilde{g}_{(i,1)} \widetilde{\rho}_{1}^{-1} \widetilde{g}_{(i,2)} \widetilde{\rho}_{(i,1)} \quad 
\text{(Definition of $\widetilde{\rho}_{(i,1)}$)}\\
&=&\widetilde{g}_{(i,2)}\widetilde{g}_{(i,1)} \widetilde{\rho}_{1} \widetilde{g}_{(i,1)} \widetilde{g}_{(i,2)} \widetilde{\rho}_{(i,1)}\widetilde{\rho}_{1}^{-1} \quad 
\text{(Commutativity of $\widetilde{\rho}_{1}$ and $\widetilde{\rho}_{(i,1)}$)}\\
&=&(\widetilde{g}_{(i,2)}\widetilde{g}_{(i,1)} \widetilde{\rho}_{1} \widetilde{g}_{(i,1)} \widetilde{g}_{(i,2)} \widetilde{g}_{(i,1)} \widetilde{\rho}_{1}  \widetilde{g}_{(i,1)} \widetilde{\rho}_{1}^{-1} )\widetilde{\rho}_{1}^{-1} \quad \text{(Definition of $\widetilde{\rho}_{(i,1)}$)}
\\ 
&=&(\widetilde{g}_{(i,2)}\widetilde{g}_{(i,1)} \widetilde{\rho}_{1} \widetilde{g}_{(i,2)} \widetilde{g}_{(i,1)} \widetilde{g}_{(i,2)} \widetilde{\rho}_{1}  \widetilde{g}_{(i,1)})\widetilde{\rho}_{1}^{-2} \quad \text{(By the relation \eqref{A1.1})}
\\
&=&(\widetilde{g}_{(i,2)}\widetilde{g}_{(i,1)}\widetilde{g}_{(i,2)} \widetilde{\rho}_{1} \widetilde{g}_{(i,1)} \widetilde{\rho}_{1} \widetilde{g}_{(i,2)} \widetilde{g}_{(i,1)})\widetilde{\rho}_{1}^{-2} \\
&=&(\widetilde{g}_{(i,1)}\widetilde{g}_{(i,2)}\widetilde{g}_{(i,1)} \widetilde{\rho}_{1} \widetilde{g}_{(i,1)} \widetilde{\rho}_{1} \widetilde{g}_{(i,2)} \widetilde{g}_{(i,1)})\widetilde{\rho}_{1}^{-2} \quad \text{(By the relation \eqref{A1.1})}
\\
&=& (\widetilde{g}_{(i,1)} \widetilde{g}_{(i,2)} \widetilde{\rho}_{1}\widetilde{g}_{(i,1)}\widetilde{\rho}_{1} \widetilde{g}_{(i,1)} \widetilde{g}_{(i,2)} \widetilde{g}_{(i,1)}) \widetilde{\rho}_{1}^{-2} \quad \text{(Commutativity of $\widetilde{\rho}_{1}$ and $\widetilde{\rho}_{(i,1)}$)}
\\
&=& (\widetilde{g}_{(i,1)} \widetilde{g}_{(i,2)} \widetilde{\rho}_{1}\widetilde{g}_{(i,1)}\widetilde{\rho}_{1} \widetilde{g}_{(i,2)} \widetilde{g}_{(i,1)} \widetilde{g}_{(i,2)}) \widetilde{\rho}_{1}^{-2} \quad \text{(By the relation \eqref{A1.1})}
\\
&=& (\widetilde{g}_{(i,1)} \widetilde{\rho}_{1} \widetilde{g}_{(i,2)}  \widetilde{g}_{(i,1)}\widetilde{g}_{(i,2)} \widetilde{\rho}_{1} \widetilde{g}_{(i,1)} \widetilde{g}_{(i,2)}) \widetilde{\rho}_{1}^{-2}\\
&=& (\widetilde{g}_{(i,1)} \widetilde{\rho}_{1} \widetilde{g}_{(i,1)}  \widetilde{g}_{(i,2)}\widetilde{g}_{(i,1)} \widetilde{\rho}_{1} \widetilde{g}_{(i,1)} \widetilde{g}_{(i,2)}) \widetilde{\rho}_{1}^{-2} \quad \text{(By the relation \eqref{A1.1})}
\\
&=& (\widetilde{g}_{(i,1)} \widetilde{\rho}_{1} \widetilde{g}_{(i,1)}  \widetilde{g}_{(i,2)}\widetilde{g}_{(i,1)} \widetilde{\rho}_{1} \widetilde{g}_{(i,1)} \widetilde{\rho}_{1}^{-1} \widetilde{g}_{(i,2)}) \widetilde{\rho}_{1}^{-1}\\
&=& (\widetilde{g}_{(i,1)} \widetilde{\rho}_{1} \widetilde{g}_{(i,1)}  \widetilde{g}_{(i,2)} \widetilde{\rho}_{(i,1)} \widetilde{g}_{(i,2)}) \widetilde{\rho}_{1}^{-1}
\quad \text{(Definition of $\widetilde{\rho}_{(i,1)}$)}
\\
&=& \widetilde{g}_{(i,1)} \widetilde{\rho}_{1} \widetilde{g}_{(i,1)} \widetilde{\rho}_{1}^{-1} \widetilde{g}_{(i,2)} \widetilde{\rho}_{(i,1)} \widetilde{g}_{(i,2)} \quad 
\text{(Commutativity of $\widetilde{\rho}_{1}$ and $\widetilde{\rho}_{(i,1)}$)}\\
&=& \widetilde{\rho}_{(i,1)} \widetilde{g}_{(i,2)} \widetilde{\rho}_{(i,1)} \widetilde{g}_{(i,2)}\quad \text{(Definition of $\widetilde{\rho}_{(i,1)}$)}.
\end{eqnarray*}

Then we have 
\[
\widetilde{g}_{(i,2)}\widetilde{\rho}_{(i,1)}\widetilde{g}_{(i,2)}\widetilde{\rho}_{(i,1)} = \widetilde{\rho}_{(i,1)}\widetilde{g}_{(i,2)}\widetilde{\rho}_{(i,1)}\widetilde{g}_{(i,2)}
\]
By the above relation and the definition of $\widetilde{\rho}_{(i,2)}$, we have 
\[
\widetilde{\rho}_{(i,2)}\widetilde{\rho}_{(i,1)} = \widetilde{g}_{(i,2)}\widetilde{\rho}_{(i,1)}\widetilde{g}_{(i,2)} = \widetilde{\rho}_{(i,1)}\widetilde{g}_{(i,2)}\widetilde{\rho}_{(i,1)}\widetilde{g}_{(i,2)}\widetilde{\rho}_{(i,1)}^{-1} = \widetilde{\rho}_{(i,1)}\widetilde{\rho}_{(i,2)}.
\]
and $\widetilde{\rho}_{(i,j)}\widetilde{\rho}_{(i,j+1)}=\widetilde{\rho}_{(i,j+1)}\widetilde{\rho}_{(i,j)}$.

By the definition of $\widetilde{\rho}_{(i,j)}$, it is written by products of 
$\widetilde{g}_{1}, \widetilde{g}_{1^{*}}, \widetilde{g}_{(i,k)}$ for $k=1,\dots j$. 
Then $\widetilde{\rho}_{(i,j)}$ commutes with $\widetilde{g}_{(i,j+2)}$ by \eqref{W1.0}.
Then we have $\widetilde{\rho}_{(i,j+2)}\widetilde{\rho}_{(i,j)}=\widetilde{\rho}_{(i,j)}\widetilde{\rho}_{(i,j+2)}$. 
Therefore, by same calculations,
we obtain the commutativity of $\widetilde{\rho}_{(i,j)}$ and $\widetilde{\rho}_{(i,k)}$ for 
$j,k=1,2,\cdots, a_{i}-1$.

Similar calculations show the elements $\widetilde{\rho}_{v}, v\in T_{A}$ satisfy the relation \eqref{E c}.
\qed
\end{pf}
To construct the inverse homomorphism, let us define the
element $g_{1^{*}}$ of the group $G(\widetilde{R}_A)$ by
\begin{equation}
g_{1^{*}}:=g_{1}^{-1}\rho_{1}.
\end{equation}

\begin{lem}
The elements $g_{v}$, $v\in T_A$ and $g_{1^{*}}$ of $G(\widetilde{R}_A)$
satisfy the relations in Definition \ref{def G'}. 
\end{lem}

\begin{pf}
First, we show the relation \eqref{A1.0}.
We only have to check that $g_{1^{*}}g_{(i,j)}=g_{(i,j)}g_{1^{*}}$ for $j=2,\dots, a_{i}-1$.
By the relation \eqref{E 3}, we have 
\begin{equation}
g_{1^{*}}g_{(i,j)}=g^{-1}_{1}\rho_{1}w_{(i,j)}=g^{-1}_{1}g_{(i,j)}\rho_{1}=g_{(i,j)}g^{-1}_{1}\rho_{1}
=g_{(i,j)}g_{1^{*}}.
\end{equation}

Second, we show the relation \eqref{A1.1}.
We only have to show the relation
\begin{equation*}
g_{(i,1)} g_{1^{*}} g_{(i,1)} = g_{1^{*}} g_{(i,1)} g_{1^{*}},
\end{equation*}
which follows from the calculation
\begin{eqnarray*}
g_{1^{*}} g_{(i,1)} g_{1^{*}}
&=& g_{1}g_{(i,1)}(g_{(i,1)} g_{1} g_{1^{*}} g_{(i,1)})g_{1^{*}}\\
&=& g_{1}g_{(i,1)}(g_{(i,1)} \rho_{1} g_{(i,1)} \rho_{1}^{-1} \rho_{1})g_{1^{*}}\\
&=& g_{1}g_{(i,1)}(\rho_{(i,1)} \rho_{1} )g_{1^{*}}\\
&=&g_{1}g_{(i,1)}(g_{1} \rho_{(i,1)} g_{1})g_{1^{*}}\\
&=& g_{1}g_{(i,1)}(g_{1} g_{(i,1)} \rho_{1} g_{(i,1)} \rho_{1}^{-1} g_{1})g_{1^{*}} \\
&=& g_{1}g_{(i,1)}(g_{1} g_{(i,1)} g_{1} g_{1^{*}} g_{(i,1)} (g_{1^{*}})^{-1} g_{1}^{-1} g_{1})g_{1^{*}} \\
&=& g_{1}g_{(i,1)}(g_{(i,1)} g_{1} g_{(i,1)} g_{1^{*}} g_{(i,1)} (g_{1^{*}})^{-1})g_{1^{*}}\\
&=& g_{(i,1)} g_{1^{*}}g_{(i,1)}.
\end{eqnarray*}

Finally, we show the relation \eqref{A2}.
Since the elements $g_{(i,1)}$ and $\rho_{1}$ satisfy the relation
\begin{equation}
g_{(i,1)} \rho_{1} g_{(i,1)} = \rho_{(i,1)} \rho_{1},
\end{equation}
by commutativity of $\rho_{(i,1)}$ and $\rho_{1}$, we have
\begin{equation*}
g_{(i,1)}\rho_{1}g_{(i,1)}\rho_{1}=\rho_{(i,1)} \rho_{1} \rho_{1} =\rho_{1} \rho_{(i,1)} \rho_{1} = \rho_{1} g_{(i,1)} \rho_{1} g_{(i,1)}. 
\end{equation*}
The relation \eqref{A3} is clearly satisfied since it is a part of the relation~\eqref{E 3}.
\qed
\end{pf} 
Therefore we have finished the proof of Theorem \ref{G isom G'}. 
\qed
\end{pf}
There exists the following commutative diagram of groups
\begin{equation}
\begin{CD}
G(\widetilde{T}_A)@>>> G(\widetilde{R}_A)\\
@VV{}V @VV{}V \\
W(\widetilde{T}_A)@>>> W(R_A)\ltimes K_0(R_A)
\end{CD},
\end{equation}
where the upper horizontal homomorphism is the isomorphisms in Theorem \ref{main},
the lower horizontal homomorphism is the isomorphisms in Theorem \ref{G isom G'},
the left vertical homomorphisms is the one in Proposition \ref{W' isom G' bar}
and finally, the right vertical homomorphism is the one in Corollary \ref{surj:5.8}. 
\section{Autoequivalence group}
In this section, we compare the cuspidal Artin group $G(\widetilde{T}_A)$ with a subgroup of autoequivalence group
for the derived category of the $2$-Calabi--Yau completion of $k\widetilde{\TT}_{A,\Lambda}$ generated by some 
spherical twist functors.

\begin{defn}\label{def:6.1}
Put $\A:=k\widetilde{\TT}_{A,\Lambda}$ and consider it as a dg $k$-algebra concentrated in the degree $0$. 
Let $\Theta_{\A}$ be the cofibrant 
replacement of the complex $\RR {\rm Hom}_{\A\otimes_{k} \A^{op}}(\A, \A\otimes_{k} \A^{op})$. The {\it $2$-Calabi--Yau 
completion} (or {\it derived $2$-preprojective algebra}) of $\A$ is the following tensor dg $k$-algebra:
\begin{equation}
\displaystyle
\Pi_{2}(\A):=\A\bigoplus_{n\in \NN} \underbrace{(\Theta_{\A}[1] \otimes_{\A}\cdots\otimes_{\A} \Theta_{\A}[1])}_{n-times}.
\end{equation}
\end{defn}
\begin{rem}
Since $k\widetilde{\TT}_{A,\Lambda}$ is a directed finite dimensional algebra over the {\it field} $k$ of global dimension two, 
the above definition agrees with the original one in \cite{keller}.
\end{rem}
Let $\D(\Pi_2(\A))$ be the derived category of dg $\Pi_2(\A)$-modules.
Note that we have a natural functor $\D(k\widetilde{\TT}_{A,\Lambda})\longrightarrow \D(\Pi_2(\A))$ given by the restriction along the projection onto the first component $\Pi_2(\A)\longrightarrow \A =k\widetilde{\TT}_{A,\Lambda}$.
Therefore we shall often regard $M\in \D(k\widetilde{\TT}_{A,\Lambda})$ also as a dg $\Pi_2(\A)$-module.
Let $\check{\D}_{A,\Lambda}$ be the smallest full triangulated subcategory of $\D(\Pi_2(\A))$ 
containing $k\widetilde{\TT}_{A,\Lambda}$, closed under isomorphisms and taking direct summand.
\begin{prop}\label{prop:K-isom}
The functor $\D(k\widetilde{\TT}_{A,\Lambda})\longrightarrow \D(\Pi_2(\A))$ induces 
an isomorphism of abelian groups $K_0(\widetilde{R}_A)=K_0(\D^b(k\widetilde{\TT}_{A,\Lambda}))\cong K_0(\check{\D}_{A,\Lambda})$.
\end{prop}
\begin{pf}
The statement follows from the fact that the triangulated category $\D^b(k\widetilde{\TT}_{A,\Lambda})$ is equivalent to
the smallest full triangulated subcategory of  $\D(k\widetilde{\TT}_{A,\Lambda})$ 
containing $k\widetilde{\TT}_{A,\Lambda}$, closed under isomorphisms and taking direct summand.
\qed
\end{pf}
\begin{prop}\label{prop:63}
For any $X,Y\in \D^b(k\widetilde{\TT}_{A,\Lambda})$, there is a canonical isomorphism in $\D^b(k):$
\begin{equation}
\RR{\rm Hom}_{\check{\D}_{A,\Lambda}}(X,Y)\cong \RR{\rm Hom}_{\D^b(k\widetilde{\TT}_{A,\Lambda})}(X,Y)\oplus \RR{\rm Hom}_{\D^b(k\widetilde{\TT}_{A,\Lambda})}(Y,X)^*[-2].
\end{equation}
\end{prop}
\begin{pf}
This is a direct consequence of Lemma~4.4 b) in \cite{keller}.
\qed
\end{pf}
\begin{cor}\label{cor:65}
Under the isomorphism $K_0(\widetilde{R}_A)\cong K_0(\check{\D}_{A,\Lambda})$ in Proposition~\ref{prop:K-isom}, 
the Euler form $\chi_{\check{\D}_{A,\Lambda}}$ is identified with the Cartan form $I_{\D^b(k\widetilde{\TT}_{A,\Lambda})}$.
\end{cor}
\begin{pf} 
It follows from Proposition~\ref{prop:63} that we have
\[
\chi_{\check{\D}_{A,\Lambda}}([X],[Y])=\chi_{\D^b(k\widetilde{\TT}_{A,\Lambda})}([X],[Y])+\chi_{\D^b(k\widetilde{\TT}_{A,\Lambda})}([Y],[X])=I_{\D^b(k\widetilde{\TT}_{A,\Lambda})}([X],[Y])
\]
for any $X,Y\in \D^b(k\widetilde{\TT}_{A,\Lambda})$ regarded as objects in 
$\check{\D}_{A,\Lambda}$.
\qed
\end{pf}
Recall the definitions of spherical objects and spherical twist functors and their properties in Seidel--Thomas \cite{Seidel-Thomas}.
\begin{defn}\label{def:6.2}
An object $S\in \check{\D}_{A,\Lambda}$ is called a {\it $2$-spherical object} if the following conditions are satisfied:
\begin{enumerate}
\item There exists an isomorphism in $\D^{b}(k)$:
\begin{equation}
\RR{\rm Hom}_{\check{\D}_{A,\Lambda}}(S,S)\cong
k\oplus k[-2]
\end{equation}
\item For all $X\in \check{\D}_{A,\Lambda}$, the composition induces the following perfect pairing:
\begin{equation}\label{def:6.2(6.3)}
{\rm Hom}_{\check{\D}_{A,\Lambda}}(X,S[2])\otimes_{k} {\rm Hom}_{\check{\D}_{A,\Lambda}}(S,X)
\rightarrow {\rm Hom}_{\check{\D}_{A,\Lambda}}(S,S[2]) \cong k.
\end{equation}
\end{enumerate}
\end{defn}
\begin{defn}
Let $S$ be a spherical object in $\check{\D}_{A,\Lambda}$ and $X$ a object in $\check{\D}_{A,\Lambda}$.
Define $T_S X\in\check{\D}_{A,\Lambda}$ by the cone of the evaluation morphism $ev$
\begin{equation}
\RR{\rm Hom}_{\check{\D}_{A,\Lambda}}(S,X)\otimes^\LL S\stackrel{ev}{\longrightarrow} X.
\end{equation}
Similarly, define $T^-_S X\in\check{\D}_{A,\Lambda}$ by the $-1$-translation of the cone of the evaluation morphism $ev^*$
\begin{equation}
X\stackrel{ev^*}{\longrightarrow} \RR{\rm Hom}_{\check{\D}_{A,\Lambda}}(X,S)^*\otimes^\LL S.
\end{equation}
The operations $T_{S}$ and $T^{-}_{S}$ define endo-functors on $\check{\D}_{A,\Lambda}$, which are
called the {\it spherical twist} functors.
\end{defn}
We collect some basic properties of the spherical twist functors.
In particular, it turns out that the spherical twist functors are autoequivalences on $\check{\D}_{A,\Lambda}$.
\begin{prop}
Let $S$ be a spherical object in $\check{\D}_{A,\Lambda}$.
\begin{enumerate}
\item
For an integer $i\in\ZZ$, we have $T_{S[i]}\cong T_S$.
\item 
We have $T^-_S T_S\cong {\rm Id}_{\check{\D}_{A,\Lambda}}$ and $T_S T^-_S\cong {\rm Id}_{\check{\D}_{A,\Lambda}}$.
\item
We have $T_{S}S\cong S[-1]$.
\item
For any spherical object $S'$, we have
\begin{equation}\label{64}
T_{S}T_{S'}\cong T_{T_{S}S'}T_{S}.
\end{equation}
\item
For any spherical objects $S'$ such that $\RR{\rm Hom}_{\check{\D}_{A,\Lambda}}(S',S)\cong k[-1]$ in $\D(k)$, 
we have an isomorphism
\begin{equation}\label{65}
T_{S}T_{S'}S\cong S'\quad \text{in}\quad \check{\D}_{A,\Lambda}.
\end{equation}
\end{enumerate}
\end{prop}
\begin{pf}
The relations {\rm (i), (ii), (iii)}, the equations \eqref{64} and \eqref{65} follow from Proposition 2.10, Lemma 2.11 and Proposition 2.13 in \cite{Seidel-Thomas}. 
\qed
\end{pf}
Recall that $S_v$ is the simple $k\widetilde{\TT}_{A,\Lambda}$-module corresponding to the vertex $v\in \widetilde{T}_A$ (see Definition~\ref{def:226}), which we regard as a dg $\Pi_2(k\widetilde{\TT}_{A,\Lambda})$-module. 
\begin{prop}
For any $v\in \widetilde{T}_A$, $S_v$ is a spherical object in $\check{\D}_{A,\Lambda}$.
\end{prop}
\begin{pf}
Since $S_v$ is an exceptional object in $\D^b(k\widetilde{\TT}_{A,\Lambda})$, it follows from Propositon~\ref{prop:63}.
\qed
\end{pf}
\begin{prop}
Under the isomorphism $K_0(\check{\D}_{A,\Lambda})\cong K_0(\widetilde{R}_A)$ in Proposition~\ref{prop:K-isom}, 
the automorphism of $K_{0}(\widetilde{R}_{A})$ induced by $T_{S_{v}}$ is identified with the simple reflection 
$\widetilde{r}_{v}\in W(\widetilde{R}_A)$.
\end{prop}
\begin{pf}
By the definition of the spherical twist $T_{S_{v}}$ and Corollary~\ref{cor:65}, we have
\begin{eqnarray*}
[T_{S_{v}}X]&=&[X]-[\RR{\rm Hom}_{\check{\D}_{A,\Lambda}}(S_{v},X)\otimes^\LL S_{v}]\\
&=&[X]-\chi_{\check{\D}_{A,\Lambda}}([S_{v}],[X])[S_{v}]\\
&=&[X]-I_{\D^b(k\widetilde{\TT}_{A,\Lambda})}([S_{v}],[X])[S_{v}]=\widetilde{r}_{v}([X])
\end{eqnarray*}
for any $X\in \check{\D}_{A,\Lambda}$.
\qed
\end{pf}
\begin{defn}
Denote by ${\rm Br}(\check{\D}_{A,\Lambda})$ the subgroup of ${\rm Auteq}(\check{\D}_{A,\Lambda})$ generated by the elements
$T_{S_v}$ for $v\in\widetilde{T}_A$.
\end{defn}
\begin{thm}\label{isom G and Br}
The correspondence $\widetilde{g}_v\mapsto T_{S_{v}}$ for $v\in\widetilde{T}_A$ induces a surjective group homomorphism 
\begin{equation}
G(\widetilde{T}_A)\twoheadrightarrow {\rm Br}(\check{\D}_{A,\Lambda}). 
\end{equation}
\end{thm}
\begin{pf}
In order to simplify the notation, set $T_v:=T_{S_v}$ for all $v\in\widetilde{T}_A$.
We only need to check that the elements $T_v$ for $v\in\widetilde{T}_A$ satisfy 
the relations \eqref{A2} and \eqref{A3} since the relations \eqref{A1.0} and \eqref{A1.1} are already known by Seidel--Thomas
(Theorem~2.17 in \cite{Seidel-Thomas}).
We first give some facts useful in later discussion.
\begin{lem}\label{lem:610}
We have
\begin{subequations}
\begin{equation}
\RR{\rm Hom}_{\check{\D}_{A,\Lambda}}(S_{(i,1)},S_{(j,1)})\simeq 0,\quad \text{for all}\quad 1\le i<j\le r,
\end{equation}
\begin{equation}
\RR{\rm Hom}_{\check{\D}_{A,\Lambda}}(S_{1^*}, S_{(i,1)})\cong k[-1],\quad 
\RR{\rm Hom}_{\check{\D}_{A,\Lambda}}(S_{(i,1)},S_{1})\cong k[-1],
\end{equation}
\begin{equation}
\RR{\rm Hom}_{\check{\D}_{A,\Lambda}}(S_{1^*}, S_{1})\cong k^{\oplus 2}[-2].
\end{equation}
In particular, the composition map
\begin{equation}
{\rm Hom}_{\check{\D}_{A,\Lambda}}(S_{(i,1)}[1],S_{1}[2])\otimes_k 
{\rm Hom}_{\check{\D}_{A,\Lambda}}(S_{1^*}, S_{(i,1)}[1]) 
\longrightarrow {\rm Hom}_{\check{\D}_{A,\Lambda}}(S_{1^*}, S_{1}[2])
\end{equation}
is non-zero.
\end{subequations}
\end{lem}
\begin{pf}
It follows from a direct calculation in $\D^b(k\widetilde{\TT}_{A,\Lambda})$ that 
\begin{equation*}
\RR{\rm Hom}_{\D^b(k\widetilde{\TT}_{A,\Lambda})}(S_{(i,1)},S_{(j,1)})\simeq 0,\quad \text{for all}\quad 1\le i<j\le r,
\end{equation*}
\begin{equation*}
\RR{\rm Hom}_{\D^b(k\widetilde{\TT}_{A,\Lambda})}(S_{1^*}, S_{(i,1)})\cong k[-1],\quad 
\RR{\rm Hom}_{\D^b(k\widetilde{\TT}_{A,\Lambda})}(S_{(i,1)},S_{1^*})\cong 0,
\end{equation*}
\begin{equation*}
\RR{\rm Hom}_{\D^b(k\widetilde{\TT}_{A,\Lambda})}(S_{(i,1)},S_{1})\cong k[-1], \quad
\RR{\rm Hom}_{\D^b(k\widetilde{\TT}_{A,\Lambda})}(S_{1},S_{(i,1)})\cong 0,
\end{equation*}
\begin{equation*}
\RR{\rm Hom}_{\D^b(k\widetilde{\TT}_{A,\Lambda})}(S_{1^*}, S_{1})\cong k^{\oplus 2}[-2],\quad 
\RR{\rm Hom}_{\D^b(k\widetilde{\TT}_{A,\Lambda})}(S_{1}, S_{1^*})\cong 0.
\end{equation*}
In particular, the composition map
\[
{\rm Hom}_{\D^b(k\widetilde{\TT}_{A,\Lambda})}(S_{(i,1)}[1],S_{1}[2])\otimes_k 
{\rm Hom}_{\D^b(k\widetilde{\TT}_{A,\Lambda})}(S_{1^*}, S_{(i,1)}[1]) 
\longrightarrow {\rm Hom}_{\D^b(k\widetilde{\TT}_{A,\Lambda})}(S_{1^*}, S_{1}[2])
\]
is non-zero. Proposition~\ref{prop:63} yields the statement. 
\qed
\end{pf}
For the later use,
we choose morphisms $\varphi^{*}_{(l,1)}$ and $\varphi_{(l,1)}$ for $l=1,\dots r$ such that 
\begin{subequations}
\begin{equation}\label{mor:6.11}
{\rm Hom}_{\check{\D}_{A,\Lambda}}(S_{1^{*}}, S_{(j,1)}[1])=k\cdot \varphi^{*}_{(l,1)},
\end{equation}
\begin{equation}\label{mor:6.12}
{\rm Hom}_{\check{\D}_{A,\Lambda}}(S_{(j,1)}[1], S_{1}[2])=k\cdot \varphi_{(l,1)}.
\end{equation}
\end{subequations}
We check the relation \eqref{A2}.
By the equation \eqref{64}, we have
\[
T_1T_{1^*}T_{(i,1)}T_1T_{1^*}T_{(i,1)}\cong T_{T_1T_{1^*}T_{(i,1)}T_1T_{1^*}S_{(i,1)}}T_1T_{1^*}T_{(i,1)}T_1T_{1^*}T_{(i,1)}.
\]
By the equations \eqref{64} and \eqref{65}, we obtain
\[
T_1T_{1^*}T_{(i,1)}T_1T_{1^*}S_{(i,1)}\cong T_1T_{1^*}T_{T_{(i,1)}S_1}T_{(i,1)}T_{1^*}S_{(i,1)} \cong T_1T_{1^*}T_{T_{(i,1)}S_1}S_{1^*}.
\]
\begin{lem}
There are the following isomorphisms in $\D^b(k):$
\begin{subequations}
\begin{equation}\label{6.10a}
\RR{\rm Hom}_{\check{\D}_{A,\Lambda}}(S_{1^*},T_{(i,1)}S_1)\cong k[-2], 
\end{equation}
\begin{equation}\label{6.10b}
\RR{\rm Hom}_{\check{\D}_{A,\Lambda}}(T_{(i,1)}S_1,S_{1^*})\cong k.
\end{equation}
\end{subequations}
\end{lem}
\begin{pf}
Note that the quasi-isomorphism $\RR{\rm Hom}_{\check{\D}_{A, \Lambda}}(S_{(i,1)}, S_{1})\simeq k[-1]$ follows from 
the quasi-isomorphism $\RR{\rm Hom}_{\check{\D}_{A, \Lambda}}(S_{1},S_{(i,1)})\simeq k[-1]$
and the isomorphism, induced by the perfect pairing in Definition \ref{def:6.2}, 
$\RR{\rm Hom}_{\check{\D}_{A, \Lambda}}(S_{1},S_{(i,1)}[2])^{*}\cong 
\RR{\rm Hom}_{\check{\D}_{A, \Lambda}}(S_{(i,1)}, S_{1})$.

By the definition of $T_{(i,1)}S_{1}$, we have the following long exact sequence:
{\small
\[
\cdots\longrightarrow
{\rm Hom}_{\check{\D}_{A,\Lambda}}(S_{1^*},S_{(i,1)}[-1])\longrightarrow 
{\rm Hom}_{\check{\D}_{A,\Lambda}}(S_{1^*},S_1)\longrightarrow 
{\rm Hom}_{\check{\D}_{A,\Lambda}}(S_{1^*},T_{(i,1)}S_1)\longrightarrow \cdots
\]
}
Since the map ${\rm Hom}_{\check{\D}_{A,\Lambda}}(S_{1^*},S_{(i,1)}[1])\longrightarrow 
{\rm Hom}_{\check{\D}_{A,\Lambda}}(S_{1^*},S_1[2])$ is injective by Lemma~\ref{lem:610},
we obtain the first isomorphism \eqref{6.10a}
The first isomorphism implies the second one \eqref{6.10b} by the isomorphism
$\RR{\rm Hom}_{\check{\D}_{A,\Lambda}}(S_{1^*},T_{(i,1)}S_1[2])^{*}\cong \RR{\rm Hom}_{\check{\D}_{A,\Lambda}}(T_{(i,1)}S_1,S_{1^*})$
induced by the perfect pairing in Definition \ref{def:6.2}.
\qed
\end{pf}
By this lemma and the equation \eqref{65}, we get
\[
T_1T_{1^*}T_{T_{(i,1)}S_1}S_{1^*}\cong T_1T_{(i,1)}S_1[1]\cong S_{(i,1)}[1].
\]
Therefore, $T_{T_1T_{1^*}T_{(i,1)}T_1T_{1^*}S_{(i,1)}}\cong T_{(i,1)}$, which gives the relation \eqref{A2}, namely, 
\begin{equation}
T_{(i,1)}T_1T_{1^*}T_{(i,1)}T_1T_{1^*}\cong T_1T_{1^*}T_{(i,1)}T_1T_{1^*}T_{(i,1)}.
\end{equation}
Next we show the relation \eqref{A3}.
For $1\le i<j\le r$, by the equation \eqref{64}, we have
\[
T_{(j,1)}T_1T_{1^*}T_{(j,1)}T^-_{1^*}T^-_{1}T_{(i,1)}\cong T_{T_{(j,1)}T_1T_{1^*}T_{(j,1)}T^-_{1^*}T^-_{1}S_{(i,1)}}T_{(j,1)}T_1T_{1^*}T_{(j,1)}T^-_{1^*}T^-_{1},
\]
and 
\begin{eqnarray*}
T_{(j,1)}T_1T_{1^*}T_{(j,1)}T^-_{1^*}T^-_{1}S_{(i,1)}&\cong& T_{(j,1)}T_1T_{T_{1^*}S_{(j,1)}}T_{1^*}T^-_{1^*}T^-_{1}S_{(i,1)}\\
&\cong& T_{(j,1)}T_1T_{T_{1^*}S_{(j,1)}}T^-_{1}S_{(i,1)}\\
&\cong& T_{(j,1)}T_{T_1T_{1^*}S_{(j,1)}}T_1T^-_{1}S_{(i,1)}\\
&\cong& T_{(j,1)}T_{T_1T_{1^*}S_{(j,1)}}S_{(i,1)}.
\end{eqnarray*}
\begin{lem}
For $1\le i<j\le r$, there are the following isomorphisms in $\D(k):$
\begin{subequations}
\begin{equation}\label{eq:6.12a}
\RR{\rm Hom}_{\check{\D}_{A,\Lambda}}(S_{(i,1)},T_1T_{1^*}S_{(j,1)})\cong 0, 
\end{equation}
\begin{equation}\label{eq:6.12b}
\RR{\rm Hom}_{\check{\D}_{A,\Lambda}}(T_1T_{1^*}S_{(j,1)},S_{(i,1)})\cong 0.
\end{equation}
\end{subequations}
\end{lem}
\begin{pf}
Note that the quasi-isomorphism 
$\RR{\rm Hom}_{\check{\D}_{A, \Lambda}}(S_{1^{*}}, S_{(j,1)})\simeq k[-1]$ follows from 
the quasi-isomorphism 
$\RR{\rm Hom}_{\check{\D}_{A, \Lambda}}(S_{(j,1)}, S_{1^{*}})\cong k[-1]$
and the isomorphism, induced by the perfect pairing in Definition \ref{def:6.2},
$\RR{\rm Hom}_{\check{\D}_{A, \Lambda}}(S_{(j,1)}, S_{1^{*}}[2])^{*}
\cong \RR{\rm Hom}_{\check{\D}_{A, \Lambda}}(S_{1^{*}}, S_{(j,1)})$.
First, by definition of $T_{1^{*}}S_{(j,1)}$ we have the following long exact sequence:
{\small
\begin{equation}\label{seq:6.14}
\cdots\longrightarrow
{\rm Hom}_{\check{\D}_{A,\Lambda}}(S_{(i,1)},S_{1^*}[-1])\longrightarrow 
{\rm Hom}_{\check{\D}_{A,\Lambda}}(S_{(i,1)},S_{(j,1)})\longrightarrow 
{\rm Hom}_{\check{\D}_{A,\Lambda}}(S_{(i,1)},T_{1^*}S_{(j,1)})\longrightarrow \cdots
\end{equation}
}
Since $\RR{\rm Hom}_{\check{\D}_{A,\Lambda}}(S_{(i,1)},S_{(j,1)})\simeq 0$, 
we have ${\rm Hom}_{\check{\D}_{A,\Lambda}}(S_{(i,1)},S_{(j,1)}[n])\cong 0$, namely,
the isomorphism ${\rm Hom}_{\check{\D}_{A,\Lambda}}(S_{(i,1)},T_{1^*}S_{(j,1)}[n])
\cong {\rm Hom}_{\check{\D}_{A,\Lambda}}(S_{(i,1)},S_{1^*}[n])$ for $n\in \ZZ$.
Therefore we have a natural isomorphism in $\D(k)$:
\begin{equation}\label{iso:6.15} 
\RR{\rm Hom}_{\check{\D}_{A,\Lambda}}(S_{(i,1)},T_{1^*}S_{(j,1)})\cong
\RR{\rm Hom}_{\check{\D}_{A,\Lambda}}(S_{(i,1)},S_{1^{*}})\simeq k[-1]. 
\end{equation}
Second, by definition of $ T_{1^{*}}S_{(j,1)}$, we have the following long exact sequence:
{\small
\begin{equation}\label{seq:6.16}
\cdots\longrightarrow
{\rm Hom}_{\check{\D}_{A,\Lambda}}(S_{1},S_{1^*}[-1])\longrightarrow 
{\rm Hom}_{\check{\D}_{A,\Lambda}}(S_{1},S_{(j,1)})\longrightarrow 
{\rm Hom}_{\check{\D}_{A,\Lambda}}(S_{1},T_{1^*}S_{(j,1)})\longrightarrow \cdots
\end{equation}
}
Since the morphism 
${\rm Hom}_{\check{\D}_{A,\Lambda}}(S_{1},S_{1^*})\longrightarrow {\rm Hom}_{\check{\D}_{A,\Lambda}}(S_{1},S_{(j,1)}[1])$ 
is surjective by Lemma~\ref{lem:610}, we obtain the following short exact sequence from the long exact sequence \eqref{seq:6.16}:
\begin{equation}\label{ex seq:6.17}
0\rightarrow {\rm Hom}_{\check{\D}_{A,\Lambda}}(S_{1},T_{1^*}S_{(j,1)})\rightarrow {\rm Hom}_{\check{\D}_{A,\Lambda}}(S_{1},S_{1^*})\rightarrow {\rm Hom}_{\check{\D}_{A,\Lambda}}(S_{1},S_{(j,1)}[1])\rightarrow 0. 
\end{equation}
Hence we also have $\RR{\rm Hom}_{\check{\D}_{A,\Lambda}}(S_{1},T_{1^*}S_{(j,1)})\cong k$. 
Finally, by definition of $T_1T_{1^*}S_{(j,1)}$, we have the following long exact sequence:
{\small
\begin{equation}
\cdots\longrightarrow
{\rm Hom}_{\check{\D}_{A,\Lambda}}(S_{(i,1)},S_{1})\longrightarrow 
{\rm Hom}_{\check{\D}_{A,\Lambda}}(S_{(i,1)},T_{1^*}S_{(j,1)})\longrightarrow 
{\rm Hom}_{\check{\D}_{A,\Lambda}}(S_{(i,1)},T_1T_{1^*}S_{(j,1)})\longrightarrow \cdots
\end{equation}
}
We want to show that ${\rm Hom}_{\check{\D}_{A,\Lambda}}(S_{(i,1)},S_{1}[1]) \cong {\rm Hom}_{\check{\D}_{A, \Lambda}}(S_{(i,1)}, T_{1^*}S_{(j,1)}[1])$.
Denote by $\psi$ the non-zero morphism in ${\rm Hom}_{\check{\D}_{A,\Lambda}}(S_{1}[1],S_{1^{*}}[1])$
which factors through $T_{1^*}S_{(j,1)}[1]$.
By the isomorphism \eqref{iso:6.15}, it is enough to check that 
$\psi_{*}:{\rm Hom}_{\check{\D}_{A,\Lambda}}(S_{(i,1)},S_{1}[1]) \longrightarrow {\rm Hom}_{\check{\D}_{A, \Lambda}}(S_{(i,1)}, S_{1^{*}}[1])$
is isomorphism where $\psi_{*}$ is the induced morphism by $\psi$.
In order to confirm this isomorphism, we only have to prove that $\psi_*$ is surjective since we have ${\rm Hom}_{\check{\D}_{A,\Lambda}}(S_{(i,1)},S_{1}[1])\cong k$ and 
${\rm Hom}_{\check{\D}_{A, \Lambda}}(S_{(i,1)}, S_{1^{*}}[1])\cong k$.
\begin{sublem}
For the morphism $\varphi_{(i,1)}$ in \eqref{mor:6.12}, we have $\psi\circ \varphi_{(i,1)}[-1]\ne 0$.
\end{sublem}
\begin{pf}
Since we have the following perfect pairing
\begin{equation}
{\rm Hom}_{\check{\D}_{A, \Lambda}}(S_{1^{*}}[1],S_{1}[3])\otimes {\rm Hom}_{\check{\D}_{A, \Lambda}}(S_{1}[1],S_{1^{*}}[1])
\longrightarrow {\rm Hom}_{\check{\D}_{A, \Lambda}}(S_{1}[1],S_{1}[3])\cong k,
\end{equation}
the $k$-module ${\rm Hom}_{\check{\D}_{A, \Lambda}}(S_{1}[1],S_{1^{*}}[1])$ is the dual of ${\rm Hom}_{\check{\D}_{A, \Lambda}}(S_{1^{*}}[1],S_{1}[3])$.
We have $\varphi^{*}_{(j,1)}[1]\circ \psi=0$, namely, $(\varphi_{(j,1)}[1]\circ \varphi^{*}_{(j,1)}[1])\circ \psi=0$ by the short exact sequence \eqref{ex seq:6.17}.  
For the morphisms $\varphi^{*}_{(l,1)}$ in \eqref{mor:6.11} and $\varphi_{(l,1)}$ in \eqref{mor:6.12},
the $k$-modules generated by $\varphi_{(l,1)}\circ \varphi^{*}_{(l,1)}$ for $l=1,\dots, r$
define different $k$-submodules of rank one in the $k$-module
${\rm Hom}_{\check{\D}_{A,\Lambda}}(S_{1^{*}}, S_{1}[2])$  of rank two.
Hence we have $\varphi_{(i,1)}[1]\circ \varphi^{*}_{(i,1)}[1]\circ \psi\ne 0$, namely, $(\varphi^{*}_{(i,1)}[1]) \circ \psi\ne 0$ for $i\ne j$.
Considering the following perfect pairing
\begin{equation*}
{\rm Hom}_{\check{\D}_{A,\Lambda}}(S_{1}[1],S_{(i,1)}[2])\otimes {\rm Hom}_{\check{\D}_{A,\Lambda}}(S_{(i,1)}, S_{1}[1])\longrightarrow
{\rm Hom}_{\check{\D}_{A,\Lambda}}(S_{(i,1)}, S_{(i,1)}[2])\cong k,
\end{equation*}  
\begin{equation*}
(\varphi^{*}_{(i,1)}[1]\circ \psi)\otimes \varphi_{(i,1)}[-1]\mapsto (\varphi^{*}_{(i,1)}[1]\circ \psi)\circ \varphi_{(i,1)}[-1],
\end{equation*}
we have $(\varphi^{*}_{(i,1)}[1]\circ \psi)\circ \varphi_{(i,1)}[-1]\ne 0$ and hence $\psi\circ \varphi_{(i,1)}[-1]\ne 0$. 
\qed
\end{pf}
Therefore we have 
$\RR{\rm Hom}_{\check{\D}_{A,\Lambda}}(S_{(i,1)},S_{1})
\cong \RR{\rm Hom}_{\check{\D}_{A,\Lambda}}(S_{(i,1)},T_{1^*}S_{(j,1)})\cong k[-1]$ 
and hence
we obtain the first isomorphism 
$\RR{\rm Hom}_{\check{\D}_{A,\Lambda}}(S_{(i,1)},T_1T_{1^*}S_{(j,1)})\cong 0$. 
We also obtain the second one \eqref{eq:6.12b}
by the isomorphism, induced by the perfect pairing in Definition \ref{def:6.2},
$\RR{\rm Hom}_{\check{\D}_{A,\Lambda}}(S_{(i,1)},T_1T_{1^*}S_{(j,1)}[2])^{*}
\cong \RR{\rm Hom}_{\check{\D}_{A,\Lambda}}(T_1T_{1^*}S_{(j,1)},S_{(i,1)})$.
\qed
\end{pf}
By this lemma, we get 
\[
T_{(j,1)}T_{T_1T_{1^*}S_{(j,1)}}S_{(i,1)}\cong T_{(j,1)}S_{(i,1)}\cong S_{(i,1)}.
\]
Therefore, we have the relation \eqref{A3}, namely, 
\[
T_{(i,1)}T_{(j,1)}T_1T_{1^*}T_{(j,1)}T^-_{1^*}T^-_{1}\cong T_{(j,1)}T_1T_{1^*}T_{(j,1)}T^-_{1^*}T^-_{1}T_{(i,1)}.
\]
We have finished the proof of the theorem.
\qed
\end{pf}
There exists the following commutative diagram of groups
\begin{equation}
\begin{CD}
G(\widetilde{T}_A)@>>>  {\rm Br}(\check{\D}_{A,\Lambda})\\
@VV{}V @VV{}V \\
W(\widetilde{T}_A)@>>> W(\widetilde{R}_A)
\end{CD},
\end{equation}
where the upper horizontal homomorphism is induced by
the above correspondence, the lower horizontal homomorphism is the composition of the morphisms 
in Proposition \ref{prop:42} {\rm (ii)} and Theorem \ref{main}, the left vertical homomorphism
is the surjective one in Theorem \ref{W' isom G' bar} and finally, the right vertical
homomorphism is induced by the correspondence 
$T_{S_{v}}\mapsto \widetilde{r}_{v}$ for $v\in \widetilde{T}_A$. 
Recall that the lower horizontal homomorphism is an isomorphism when $\chi_A\ne 0$.
We expect that if $\chi_A\ne 0$ then the upper horizontal homomorphism is also an isomorphism.

\end{document}